\documentclass[12pt]{amsart}
\usepackage{amssymb,latexsym}
\usepackage{enumerate}
\usepackage{amscd}
\usepackage{verbatim}

\theoremstyle{plain}
\newtheorem{theorem}{Theorem}
\newtheorem{corollary}[theorem]{Corollary}
\newtheorem{lemma}[theorem]{Lemma}

\theoremstyle{definition}

\theoremstyle{remark}

\numberwithin{equation}{section}
\numberwithin{theorem}{section}
\numberwithin{definition}{section}

\begin{document}

% One author
\title[Abelian Extensions]{Abelian Extensions of
Algebras\\ in Congruence-Modular Varieties
}
\author{William H. Rowan}
\address{48 Shattuck Square, \#86 \\
         Berkeley, California 94704}
\email{rowan@crl.com}
%\thanks{thanks}
% End one author

% Some input statements - I am not sure where they should
% go.
%\input mymathdefs3

\def\pair#1#2{\langle #1, #2\rangle}
\def\triple#1#2#3{\langle #1, #2, #3\rangle}
\def\triplem#1#2#3{M(#1,#2,#3)}
\def\triplemp#1{M(#1)}
\def\semitimes{\ltimes}
\def\rsemitimes{\rtimes}
\def\Mof#1{\Cal M(A)}
\def\pp#1#2{{\roman P[#1,#2]}}
\def\inclusion#1#2{{\scriptstyle[#1\hookrightarrow #2]}}

\newcommand{\N}{{\Bbb N}}
\newcommand{\Z}{{\Bbb Z}}
\newcommand{\Q}{{\Bbb Q}}
\newcommand{\R}{{\Bbb R}}
\newcommand{\A}{{\Bbb A}}
%\def\F{{\roman F}}
%\def\U{{\roman U}}

%category names
\def\Set{{\text{\bf Set}}}
\def\Mon{{\text{\bf Mon}}}
\def\Grp{{\text{\bf Grp}}}
\def\Rng{{\text{\bf Rng}}}
\def\CRng{{\text{\bf CRng}}}
\def\Ab{{\text{\bf Ab}}}
\def\Pnt{{\text{\bf Pnt}}}
\def\Alg{{\text{\bf Alg}}}
\def\Smgrp{{\text{\bf Smgrp}}}
\def\Ov{{\text{\bf Ov}}}
\def\Pre{{\text{\bf Pre}}}
\def\Lat{{\text{\bf Lat}}}
\def\LatX{{\text{\bf Lat-X}}}
\def\Modlat{{\text{\bf Modlat}}}
\def\Typ{{\text{\bf Typ}}}
\def\catClo{{\text{\bf Clo}}}
\def\Mod{{\text{\bf Mod}}}
\def\Bimod{{\text{\bf Bimod}}}
\def\Lie{{\text{\bf Lie}}}
\def\Nonass{{\text{\bf Nonass}}}
\def\Rep{{\text{\bf Rep}}}
\def\Rngd{{\text{\bf Rngd}}}
\def\Algd{{\text{\bf Algd}}}
\def\Cat{{\text{\bf Cat}}}
\def\Bimod{{\text{\bf Bimod}}}
\def\Abhp{{\text{\bf AbHp}}}
\def\Clnd{{\text{\bf Clnd}}}
\def\Sktch{{\text{\bf Sktch}}}
\def\Top{{\text{\bf Top}}}
\def\Vect{{\text{\bf Vect}}}
\def\BPSet{{\text{\bf BPSet}}}
\def\Abmod{{\text{\bf Abmod}}}
\def\Unif{{\text{\bf Unif}}}
\def\HUnif{{\text{\bf HUnif}}}
\def\CHUnif{{\text{\bf CHUnif}}}
\def\CU{{\text{\bf CU}}}
\def\HCU{{\text{\bf HCU}}}
\def\CHCU{{\text{\bf CHCU}}}

%other category-theoretic stuff
\def\natur{{\overset{\scriptscriptstyle\bullet}\to\to}}
%(was \text{\bf .} instead)

% Operator names:
\newcommand{\Cg}{\operatorname{Cg}}
\newcommand{\Con}{\operatorname{Con}}
\newcommand{\Cov}{\operatorname{Cov}}
\newcommand{\Div}{\operatorname{Div}}
\newcommand{\End}{\operatorname{End}}
\newcommand{\Eqv}{\operatorname{Eqv}}
\newcommand{\Ext}{\operatorname{Ext}}
\newcommand{\Fg}{\operatorname{Fg}}
\newcommand{\Hom}{\operatorname{Hom}}
\newcommand{\Ig}{\operatorname{Ig}}
\newcommand{\Mg}{\operatorname{Mg}}
\newcommand{\Pol}{\operatorname{Pol}}
\newcommand{\OpSemiUnif}{\operatorname{SemiUnif}}
\newcommand{\Sub}{\operatorname{Sub}}
\newcommand{\Sg}{\operatorname{Sg}}
\newcommand{\Term}{\operatorname{Term}}
\newcommand{\coker}{\operatorname{coker}}
\newcommand{\Ke}{\operatorname{Ke}}
\newcommand{\Co}{\operatorname{Co}}
\newcommand{\kerc}{\operatorname{kerc}}
\newcommand{\kerp}{\operatorname{kerp}}
\newcommand{\MyIm}{\operatorname{Im}}
\newcommand{\nat}{\operatorname{nat}}
\newcommand{\trcl}{\operatorname{tr.\ cl.}}
\newcommand{\typ}{\operatorname{typ}}
\newcommand{\ass}{\operatorname{ass}}
\newcommand{\card}{\operatorname{card}}
\newcommand{\Rab}{\operatorname{Rab}}
\newcommand{\Clo}{\operatorname{Clo}}
\newcommand{\id}{\operatorname{1}}
\newcommand{\Id}{\operatorname{1}}
\newcommand{\Idl}{\operatorname{Idl}}
\newcommand{\Fil}{\operatorname{Fil}}
\newcommand{\Der}{\operatorname{Der}}
\newcommand{\Sym}{\operatorname{Sym}}
\newcommand{\Res}{\operatorname{Res}}
\newcommand{\Ind}{\operatorname{Ind}}
\newcommand{\Pro}{\operatorname{Pro}}
\newcommand{\Ann}{\operatorname{Ann}}
\newcommand{\Ob}{\operatorname{Ob}}
\newcommand{\Ar}{\operatorname{Ar}}
\newcommand{\cod}{\operatorname{cod}}
\newcommand{\dom}{\operatorname{dom}}
\newcommand{\core}{\operatorname{core}}
\newcommand{\subcore}{\operatorname{subcore}}
\newcommand{\Fit}{\operatorname{Fit}}
\newcommand{\Frat}{\operatorname{Frat}}
\newcommand{\Lim}{\operatorname{Lim}}
\newcommand{\Colim}{\operatorname{Colim}}
\newcommand{\rR}{\operatorname{R}}
\newcommand{\rU}{\operatorname{U}}
\newcommand{\rClndD}{\operatorname{Clnd}^{\Bbb D}}
\newcommand{\rClnd}{\operatorname{Clnd}}
\newcommand{\rSktchD}{\operatorname{Sktch}^{\Bbb D}}
\newcommand{\rSktch}{\operatorname{Sktch}}
\newcommand{\Reg}{\operatorname{Reg}}
\newcommand{\DDown}{\operatorname{Ddn}}
\newcommand{\DUp}{\operatorname{Dup}}
\newcommand{\Gr}{\operatorname{Gr}}
\newcommand{\Step}{\operatorname{Step}}
\newcommand{\nullset}{\{\}}
\newcommand{\Ass}{\operatorname{Ass}}
\newcommand{\cov}{\operatorname{cov}}
\newcommand{\Aut}{\operatorname{Aut}}
\newcommand{\Nor}{\operatorname{Nor}}
\newcommand{\I}{\operatorname{I}}
\newcommand{\OpUnif}{\operatorname{Unif}}
\renewcommand{\subjclassname}{\textup{2000} Mathematics Subject
     Classification}

\def\Var{{\mathbf V}}
\def\Str{{\hbox{\bf str}}}
\def\Spl{{\hbox{\bf ret}}}
\def\Ret{{\hbox{\bf ret}}}

%special notations
\def\lsil{\lbrack\!\lbrack}
\def\rsil{\rbrack\!\rbrack}

%hyphenation stuff
\def\congruence{on\-gru\-ence\discretionary{-}{}{-}}
\def\conM/{c\congruence mod\-u\-lar}
\def\ConM/{C\congruence mod\-u\-lar}
\def\conD/{c\congruence dis\-trib\-u\-tive}
\def\ConD/{C\congruence dis\-trib\-u\-tive}
\def\conP/{c\congruence per\-mut\-a\-ble}
\def\ConP/{C\congruence per\-mut\-a\-ble}
\def\conMity/{\conM/\-i\-ty}
\def\ConMity/{\ConM/\-i\-ty}
\def\conDity/{c\congruence dis\-trib\-u\-tiv\-i\-ty}
\def\ConDity/{C\congruence dis\-trib\-u\-tiv\-i\-ty}
\def\conPity/{c\congruence per\-mut\-a\-bil\-i\-ty}
\def\ConPity/{C\congruence per\-mut\-a\-bil\-i\-ty}
\def\usprv/{un\-der\-ly\-ing-set-pre\-ser\-ving}

% some common latin abbreviations
\def\ie/{{i.e.}}
\def\Ie/{{I.e.}}
\def\eg/{{e.g.}}
\def\Eg/{{E.g.}}
\def\etc/{{etc.}}

% power set notation
\def\pow{{\Cal P}}

% arrows underneath formulas
\newdimen\mysubdimen
\newbox\mysubbox
\def\subright#1{\mathpalette\subrightx{#1}}
\def\subleft#1{\mathpalette\subleftx{#1}}
\def\subrightx#1#2{\subwhat#1{#2}\rightarrowfill}
\def\subleftx#1#2{\subwhat#1{#2}\leftarrowfill}

\def\subwhat#1#2#3{{
\setbox\mysubbox=\hbox{#3}
\mysubdimen=\wd\mysubbox
\setbox\mysubbox=\hbox{$#1#2$}
\ifnum\mysubdimen>\wd\mysubbox
\vtop{
\hbox to\mysubdimen{\hfil\box\mysubbox\hfil}
\nointerlineskip
\hbox{#3}}
\else
\mysubdimen=\wd\mysubbox
\vtop{
\box\mysubbox
\nointerlineskip
\hbox to\mysubdimen{#3}}
\fi
}}

\def\cite#1{[#1]}
\newcommand{\M}{M}
\renewcommand{\A}{A}
\newcommand{\E}{E}
\newcommand{\PI}{\pi}
\def\makemod#1{{#1^*}}
\newcommand{\factset}{\mathrm f}
\renewcommand{\vec}{\mathbf}

\date{\today}

\begin{abstract}
We define abelian extensions of algebras in
congruence-modular varieties. The theory is sufficiently
general that it includes, in a natural way, extensions of
$R$-modules for a ring $R$. We also define a cohomology
theory, which we call clone cohomology, such that the
cohomology group in dimension one is the group of
equivalence classes of extensions.
\end{abstract}
\maketitle

\section*{Introduction}

The theory of abelian extensions of algebras has a long
history, going from abelian extensions of groups and
modules to abelian extensions of algebras in arbitrary
varieties. Extensions of modules and abelian extensions of
groups and Lie algebras are standard topics in treatises
on Homological Algebra. The similarity of the situation
for groups and Lie algebras points to a common
generalization, which is our goal in this paper.
Previous attempts at such a generalization, and related
work, can be found in \cite{6}, \cite{9}, \cite{1},
\cite{3}, \cite{2}, \cite{12}, \cite{13},
\cite{5}, \cite{16}.

The most general treatment of which we are aware, in
\cite{3}, treats abelian extensions in the generality of
an algebra in an arbitrary variety of algebras, and a
Beck module over that algebra. (A \emph{Beck module} over
$A$ in a variety $\mathbf V$ is an abelian group object in
the category $(\mathbf V\downarrow A)$ of algebras over
$A$.)

Our treatment is less general because we require $\mathbf V$ to be congruence-modular. However, the apparatus of
commutator theory, in particular, the existence of a
difference term, allows us to prove an important lemma
(lemma~\ref{T:keylemma}) which lets us work with a simpler
and more conceptual definition of an extension than that
in
\cite{3}. We are also able to treat module extensions
with the same theory, in a fairly natural way.

In both cases, there is an associated cohomology theory,
the cohomology group in dimension one being the group of
equivalence classes of extensions. We have shown that these
groups are isomorphic for $\mathbf V$ congruence-modular,
but will not give the proof in this paper as it is quite
tedious. The theory in \cite{3} (and also \cite{1}
and \cite{2}) is called \emph{comonadic cohomology},
and we call ours \emph{clone cohomology}.

We frame our theory not in terms of algebras over $A$,
and Beck modules over $A$, but of objects in equivalent
categories we call the category of $A$-overalgebras and
the category of abelian group $A$-overalgebras. We will
briefly explain our reasons for doing so at the end of
\S\ref{S:categories}, where these categories are defined.

The theory of abelian extensions is interesting because
it organizes, from a certain point of view, the possible
structures of a class of algebras, related to $A$ (or an
$A$-overalgebra $Q$) and an abelian group $A$-overalgebra
$M$, into an abelian group, functorial in a way that we
will discuss in \S\ref{S:bifunctor}. Also, that abelian
group is a cohomology group for a suitable cohomology
theory derived from $A$ and $M$.

After a section of preliminaries, \S\ref{S:categories} of
the paper defines the categories of $A$-overalgebras and
abelian group $A$-overalgebras. \S\ref{S:ringoids}
sketches the theory of enveloping ringoids very briefly,
and is included to show how that theory can be applied to
constructing abelian group $A$-overalgebras free on an
$A$-tuple of sets of generators.
\S\ref{S:extensions} defines abelian extensions and
performs some preliminary analyses of them.
\S\ref{S:factorsets} defines and explores the formalism of
factor sets of extensions.
\S\ref{S:equivalence} introduces the definition of
equivalence of extensions, and defines the set $\mathbf
E_{\mathbf V}(A,M)$ of equivalence classes of extensions.
\S\ref{S:homologyobject} then shows how the set of
equivalence classes of extensions can be seen as a
cohomology group.
\S\ref{S:composition} explores composition operations
between abelian extensions and homomorphisms, and
\S\ref{S:addition} discusses the group law in the set of
equivalence classes of extensions. \S\ref{S:bifunctor}
contains a reformulation of the definition that defines
$\mathbf E_{\mathbf V}(Q,M)$ for an
$A$-overalgebra $Q$ and abelian group $A$-overalgebra $M$,
functorially in $Q$ and $M$. \S\ref{S:modules} shows how
module extensions can be treated. \S\ref{S:clone} presents
a cohomology theory we call \emph{clone cohomology}
because its definition intimately involves the clone of
the variety
$\mathbf V$ to which $A$, $Q$, and $M$ all belong.
\S\ref{S:relative} explores varying the variety $\mathbf
V$ used in defining clone cohomology, giving
\emph{relative clone cohomology}. Finally, we pose a
number of questions that seem important to ask about this
theory, and about the relationship of clone cohomology
and comonadic cohomology.

\section*{0. Preliminaries}

\subsection{Category theory}We follow
\cite{10} in terminology and notation.

\subsection{Homological algebra}We assume a familiarity with concepts and conventions of
homological algebra, such as can be found in \cite{11},
\cite{8}, and \cite{17}.

\subsection{Universal algebra}We assume the basic
definitions of universal algebra, such as can be found in
\cite{4}, are known to the reader.  Unlike some
authors, we admit the possibility that an algebra can
have an empty underlying set.

The kernel congruence of a homomorphism $f$ will be
denoted by $\ker f$. The other sort of kernel,
$f^{-1}(0)$, will be denoted by $\Ke f$.

We denote the greatest and least congruences of $A$ by
$\top_A$ and $\bot_A$, and the identity homomorphism of
$A$ by $1_A$. Basic operations or term operations of $A$
will occasionally be denoted by $\omega^A$ or $t^A$, but
we almost always drop the superscript.

If $A$ stands for an algebra, we use $U(A)$ to stand
for the underlying set of the algebra.

\subsection{Clones}A \emph{clone} is an $\mathbb N$-tuple
of sets
$V_n$ (the
$n$-ary elements of the clone $V$) such that for each $n$,
and each $i$ with $1\leq i\leq n$, there is an element
$\pi^V_{in}\in V_n$, called the \emph{$i^{\text{th}}$ of
$n$ projection}, and such that, for each $n$, each
$n'$-tuple $\mathbf v$ of elements of $V_n$, and each $v'\in
V_{n'}$, there is an element $v'\mathbf v\in V_n$, called
the \emph{clone composite} of $v'$ and $\mathbf v$,
satisfying
\begin{enumerate}
\item $\pi^V_{in}\mathbf v=v_i$,
\item $v\langle\pi^V_{1n},\ldots,\pi^V_{nn}\rangle = v$,
and
\item $u(\mathbf v\mathbf w)=(u\mathbf v)\mathbf w$,
\end{enumerate}
whenever the relevant compositions are defined.
($\mathbf v\mathbf w$ stands for $\langle v_1\mathbf w,\ldots,v_n\mathbf w\rangle$ if $\mathbf v$ is an $n$-tuple.)

As an example of a clone, given a set $S$, we have the
\emph{clone of $S$}, denoted by $\Clo S$.
$\Clo_n S$ is the set of $n$-ary functions from $S$ to
$S$,
$\pi^{\Clo S}_{in}$ is the $n$-ary function on $S$
choosing the $i^{\text{th}}$ of its $n$ arguments, and
given $n$, $n'$, an $n'$-tuple of $n$-ary functions
$\mathbf f$, and an $n'$-ary function $f'$, the clone
composite $f'\mathbf f$ is the function defined by
\[\mathbf s\mapsto f'(f_1(\mathbf s),\ldots,f_{n'}(\mathbf
s)).\]

Another example of a clone is the clone of a variety
$\mathbf V$, denoted by $\Clo\mathbf V$. Elements of
$\Clo_n\mathbf V$ are equivalence classes of $n$-ary term
operations of the algebras in $\mathbf V$, where terms $t$
and $t'$ are equivalent if $t(\mathbf x)=t'(\mathbf x)$ is
an identity of $\mathbf V$.

If $V$, $V'$ are clones, a \emph{homomorphism of clones}
from $V$ to $V'$ is an $\mathbb N$-tuple $f$ of functions
$f_n:V_n\to V'_n$, such that for all $i$ and $n$,
$f_n(\pi^V_{in})=\pi^{V'}_{in}$, and for all $n$ and
$n'$, $\mathbf v\in V^{n'}_n$, and $v'\in V_{n'}$, we have
$f_{n'}(v')f_n(\mathbf v)=f_n(v'\mathbf v)$.

An algebra $A$ in a variety $\mathbf V$ is the same as a
clone homomorphism from $\Clo V$ to $\Clo S$, where $S$
is the underlying set of $A$.

A clone $V$ can be viewed as a category with one object
for each natural number. The arrows from $n$ to $n'$ are
$n'$-tuples of elements of $V_n$, and the identity of $n$
is $\langle\pi^V_{1n},\ldots,\pi^V_{nn}\rangle$. In the
resulting category, $n$ is the $n$-fold direct power of
$1$. Thus, the category is what is often called a
\emph{theory}.

\subsection{The modular commutator}In a
congruence-modular variety $\mathbf V$ (i.e., such that for
all $A\in\mathbf V$, $\Con A$ is a modular lattice) the
congruence lattices admit a binary operation, called the
\emph{commutator}, that generalizes some well-known
operations such as the commutator of two normal subgroups
of a group. A comprehensive treatment of the theory of
this operation, and related matters, can be found in
\cite{7}.

The commutator of two congruences $\alpha$, $\beta\in\Con
A$ is denoted by $[\alpha,\beta]$. A congruence $\alpha$
is said to be \emph{abelian} if $[\alpha,\alpha]=\bot_A$.

One definition of the commutator is as follows: If
$A\in\mathbf V$, a congruence-modular variety of algebras,
and $\theta$, $\psi\in\Con A$, then $[\theta,\psi]$ is the
least congruence such that for all $n$, for all
$(n+1)$-ary terms $t$, for all $a$, $a'\in A$ such that
$a\mathrel\theta a'$, and for all $\mathbf b$, $\mathbf c\in
A^n$ such that $b_i\mathrel\psi c_i$ for all $i$, we have
\[t(a,\mathbf b)\mathrel{[\theta,\psi]}t(a,\mathbf c)\]
implies
\[t(a',\mathbf b)\mathrel{[\theta,\psi]}t(a',\mathbf c).\]

\subsection{Difference terms}If $\mathbf V$
is a congruence-modular variety of algebras, a ternary
term $d$ is called a \emph{difference term} for
$\mathbf V$ if
\begin{enumerate}
\item $d(x,x,y)=y$ is an identity of $\mathbf V$, and
\item for all $A\in\mathbf V$, $\theta\in\Con A$, and $x$,
$y\in A$ such that $x\mathrel\theta y$, we have
$d(x,y,y)\mathrel{[\theta,\theta]} x$.
\end{enumerate}
At least one such term exists for any congruence-modular
variety.

\subsection{$\mathbf V$-objects in a category}If $\mathbf C$ is a category, and $\mathbf V$ is a variety
of algebras, a $\mathbf V$-object in $\mathbf C$ is a pair
$\pair cF$, consisting of an object $c\in\mathbf C$, and a
contravariant functor $F:\mathbf C\to\mathbf V$, such that
$UF=\mathbf C(-,c)$, where $U:\mathbf V\to\Set$ is the
forgetful functor. If $\pair cF$ and $\pair{c'}{F'}$ are
$\mathbf V$-objects, a \emph{homomorphism of $\mathbf V$-objects from $\pair cF$ to $\pair{c'}{F'}$} is an arrow
$f:c\to c'$ such that for each object $d\in\mathbf C$, the
function $\mathbf C(d,f):\mathbf C(d,c)\to\mathbf C(d,c')$ is the
underlying function of a (necessarily unique, since $U$
is faithful) arrow $\bar f:F(d)\to F'(d)$.

$\mathbf V$-objects in $\mathbf C$, and the homomorphisms
between them, form a category in an obvious manner, which
we denote by $\mathbf V[\mathbf C]$.

\section{The Categories $A$-$\Set$, $\Ov[A,\mathbf V]$,
and
$\Ab[A,\mathbf V]$}\label{S:categories}

\subsection{The category of $A$-sets}Let $A$ be an algebra. We define an
\emph{$A$-set} to be a $U(A)$-tuple of sets. If $S$ is an
$A$-set, we will denote the member of the tuple
corresponding to an element $a\in A$ by $_aS$. If $S$,
$S'$ are $A$-sets, an \emph{$A$-function} from $S$ to $S'$
is a $U(A)$-tuple $f$ such that the element corresponding
to each $a\in A$, which we denote by $_af$, is a function
from $_aS$ to $_aS'$. We write $f:S\to S'$. $A$-sets and
$A$-functions form a category, $A$-$\Set$, in an obvious
manner.

\subsection{Overalgebras}If $A$ is an algebra, an
\emph{$A$-overalgebra} is an
$A$-set $Q$, such that for each $n$, and each $n$-ary basic
operation $\omega$ of the type of $A$, there
is a $U(A)^n$-tuple of functions $\omega^Q_{\vec
a}:{_{a_1}Q}\times\ldots\times{_{a_n}Q}\to{_{\omega(\vec
a)}Q}$.

In what follows, we will write $_{\vec a}Q$ for the product
${_{a_1}Q}\times\ldots\times{_{a_n}Q}$. Thus,
$\omega^Q_{\vec a}:{_{\vec a}Q}\to{_{\omega(\vec a)}Q}$.

If $Q$, $Q'$ are $A$-overalgebras, a \emph{homomorphism of
$A$-overalgebras from $Q$ to $Q'$} is an $A$-function
$f:Q\to Q'$ such that for each $n$, each $n$-ary
basic operation $\omega$, each $\vec a\in A^n$, and
each $\vec q\in{_{\vec a}Q}$, we
have
\[_{\omega(\vec a)}f(\omega^Q_{\vec a}(\vec q))
=\omega^{Q'}_{\vec a}({_{\vec a}f(\vec q)}),
\]
where $_{\vec a}f(\vec q)$ stands for
$\langle{_{a_1}f(q_1)},\ldots,{_{a_n}f(q_n)}\rangle$.

$A$-overalgebras and their homomorphisms form a category
in an obvious manner, which we denote by $\Ov[A]$.

For an example of an $A$-overalgebra, let $\pair B\pi$ be
an object of the ``comma category''
$(\Omega$-$\Alg\downarrow A)$ of algebras over $A$. That
is, let $B$ be an algebra of the same type $\Omega$ as $A$,
and let $\pi:B\to A$ be a homomorphism. Then we define the
$A$-overalgebra $\lsil B,\pi\rsil$ by $_a\lsil
B,\pi\rsil=\pi^{-1}(a)$ and $\omega^{\lsil
B,\pi\rsil}_{\vec a}(\vec b)=\omega^B(\vec b)$.

If $Q$ is an $A$-overalgebra, we define the \emph{total
algebra} of $Q$, denoted by $A\semitimes Q$, to be the
set of pairs $\{\,\pair aq:a\in A, q\in{_aQ}\,\}$,
provided with operations defined by
\[
\omega(\pair
{a_1}{q_1},\ldots,\pair{a_n}{q_n})=\pair{\omega(\vec
a)}{\omega^Q_{\vec a}(\vec q)};
\]
we say that an $A$-overalgebra $Q$ is \emph{totally in
$\mathbf V$}, where $V$ is a variety of algebras of the type
of $A$, if $A\semitimes Q\in\mathbf V$.

If $Q$ is an $A$-overalgebra, then accompanying the
total algebra $A\semitimes Q$ there is a
homomorphism $\pi_Q:A\semitimes Q\to A$, defined by
$\pair aq\mapsto a$.
It is clear that $\Ov[A,\mathbf V]$ is equivalent as a
category to the category $(\mathbf V\downarrow A)$. One leg
of an equivalence takes an $A$-overalgebra $Q$ to
$\pair{A\semitimes Q}{\pi_Q}$. The other leg takes an
algebra over $A$, $\pair B\pi$, to $\lsil B,\pi\rsil$.

There is an evident forgetful functor from $\Ov[A,\mathbf V]$ to $A$-$\Set$. To construct an $A$-overalgebra free
on an $A$-set $S$ (i.e., the value of the corresponding
left adjoint functor) form a free algebra $F$ on a
disjoint union of the $_aS$, (a free algebra in $\mathbf V$, that is, also known as a \emph{relatively free}
algebra), and map the generators to $A$ in the obvious
way, giving a homomorphism $\pi:F\to A$. The
$A$-overalgebra $\lsil F,\pi\rsil$ is then free on $S$.

\subsection{One-one and onto homomorphisms of
$A$-overalgebras}If $f$ is a homomorphism of $A$-overalgebras, then we say
that $f$ is \emph{one-one} if each $_af$ is one-one, and
we say that $f$ is \emph{onto} if each $_af$ is onto. If
$f$ is both one-one and onto, then it is an isomorphism
in the category of $A$-overalgebras.

\subsection{$A$-operations}Let $A$ be a set, and $\omega$ an $n$-ary operation on
$A$. If $S$ is an $A$-set, then an \emph{$A$-operation
on $S$, over $\omega$}, is a $U(A)^n$-tuple $\omega'$ of
functions $\omega'_{\mathbf a}:{_{\mathbf a}S}\to{_{\omega(\mathbf a)}S}$. We can now rephrase our
definition of an $A$-overalgebra: it is an $A$-set $Q$,
together with an $A$-operation $\omega^Q$ over
$\omega^A$ for every basic operation $\omega$ on $A$.

If $Q$ is an $A$-overalgebra totally in $\mathbf V$, and
$v\in\Clo_n\mathbf V$, then we can define an
$A$-operation $v^Q$ over $v^A$ by letting $v^Q_{\mathbf a}(\mathbf q)$ be the second component of the pair
\[
v^{A\semitimes
Q}(\pair{a_1}{q_1},\ldots,\pair{a_n}{q_n});
\]
we will thus use $v^Q$ to denote that $A$-operation,
regardless of whether $v$ is an $n$-ary basic operation,
$n$-ary term operation, or $n$-ary element of $\Clo\mathbf V$.

If $A$ is a set, and $S$ an $A$-set, then we can define
$\Clo^AS$, the \emph{clone of $A$-operations on $S$} to
be, for each $n$, the set of pairs $\pair\omega{\omega'}$
where $\omega$ is an $n$-ary operation on $A$, and
$\omega'$ is an $n$-ary $A$-operation on $S$ over
$\omega$. Projections and clone composition are easy to
define, and there is a clone homomorphism
$\pi:\Clo^AS\to\Clo A$ given by taking the first
component of each pair.

If $A$ is an algebra in a variety
$\mathbf V$, $\phi^A$ is the corresponding clone
homomorphism from $\Clo\mathbf V$ to $\Clo A$, and $Q$ is
an $A$-set, then an $A$-overalgebra structure on $Q$,
totally in $\mathbf V$, just amounts to a clone
homomorphism $\phi^Q:\Clo\mathbf V\to\Clo^AQ$, such that
$\pi\phi^Q=\phi^A$. Accordingly, to define an
$A$-overalgebra, it suffices to define $v^A$ for all $v$,
with the condition that the given definitions provide a
well-defined clone homomorphism.

\subsection{Pointed overalgebras}If $A$ is an
algebra, a \emph{pointed $A$-overalgebra} is an
$A$-overalgebra $P$, such that each $_aP$ has a basepoint
$_a*^P$, or simply $_a*$, such that for each basic
operation $\omega$, $n$-ary, and each $\mathbf a\in A^n$,
$\omega^P_{\mathbf a}({_{a_1}*},\ldots,{_{a_n}*})={_{\omega(\mathbf a)}*}$.

If $P$, $P'$ are pointed $A$-overalgebras, a \emph{
homomorphism from $P$ to $P'$} is a homomorphism of
$A$-overalgebras $f:P\to P'$ such that for each $a$,
$_af({_a*^P})={_a*^{P'}}$.

Pointed overalgebras totally in $\mathbf V$, and the
homomorphisms between them, form a category $\Pnt[A,\mathbf V]$ in an obvious manner.
As the notation suggests, it is
precisely the category $\Pnt[\Ov[A,\mathbf V]]$ of pointed
set objects in the category $\Ov[A,\mathbf V]$.

 As an example of a pointed $A$-overalgebra,
let $\alpha\in\Con A$. We define a pointed $A$-overalgebra
$\alpha^*$ by $_a\alpha^*=\{\,a':a\mathrel\alpha a'\,\}$,
$_a*=a$,
and $\omega^{\alpha^*}_{\mathbf a}(\mathbf c)=\omega(\mathbf c)$ for $\mathbf c\in{_{\mathbf a}\alpha^*}$.

If $P$ is a pointed $A$-overalgebra, then accompanying
$A\semitimes P$ and $\pi_P$ there is a homomorphism
$\iota_P:A\to A\semitimes P$, defined by
$\iota_P:a\mapsto\pair a{{_a*}}$. The triple
$\triple{A\semitimes P}{\pi_P}{\iota_P}$ can be viewed as
a commutative diagram
\[
\CD A @>\iota_P>> A\semitimes P \\
@| @VV\pi_P;V \\
A @= A\endCD
\]
in slightly different terms, it is a pointed set object
in the category of algebras over $A$.

Given such a diagram, or, a triple $\triple B\pi\iota$
with $\pi:B\to A$, $\iota:A\to B$, and $\pi\iota=1_A$, we
denote by $\lsil B,\pi,\iota\rsil$ the pointed
$A$-overalgebra with underlying $A$-overalgebra $\lsil
B,\pi\rsil$ and basepoints $\iota(a)\in{_a\lsil
B,\pi\rsil}$.

The constructions $P\mapsto\triple{A\semitimes
P}{\pi_P}{\iota_P}$ and $\triple B\pi\iota\mapsto\lsil
B,\pi,\iota\rsil$ are clearly two legs of an equivalence
between the categories $\Pnt[A,\mathbf V]$ and $\Pnt(\mathbf V\downarrow A)$.

There is an obvious forgetful functor from $\Pnt[A,\mathbf V]$ to $\Ov[A,\mathbf V]$. A corresponding free functor can
be defined as follows: given an $A$-overalgebra $Q$, form
the algebra $B=A\coprod(A\semitimes Q)$. Define $\pi:B\to
A$ by applying the universal property of the coproduct to
the homomorphisms $1_A$ and $\pi_Q$. Define $\iota:A\to
B$ as the insertion of $A$ into the coproduct. Then
$\lsil B,\pi,\iota\rsil$ is a pointed $A$-overalgebra
free on $Q$.

\subsection{Abelian group overalgebras}An \emph{abelian
group
$A$-overalgebra} is an
$A$-overalgebra $M$ such that on each $_aM$ there is the
structure of an abelian group, in such a way that the
functions $\omega^M_{\mathbf a}:{_{\mathbf a}M}\to
{_{\omega(\mathbf a)}M}$ are abelian group homomorphisms.
If $M$, $M'$ are abelian group $A$-overalgebras, a \emph{
homomorphism of abelian group $A$-overalgebras from $M$
to $M'$} is an $A$-overalgebra homomorphism $f:M\to M'$
such that each $_af$ is an abelian group homomorphism.

Abelian group $A$-overalgebras totally in $\mathbf V$, and
the homomorphisms between them, form a category in an
obvious manner, which we denote by $\Ab[A,\mathbf V]$. It
is the category of abelian group objects in $\Ov[A,\mathbf V]$.

Categorical algebraists have given the term \emph{Beck
module over $A$} to an abelian group object in the category
$(\mathbf V\downarrow A)$. Clearly, $\Ab[A,\mathbf V]$ is
equivalent to the category of Beck modules over $A$.

There is an obvious forgetful functor from $\Ab[A,\mathbf V]$ to $\Pnt[A,\mathbf V]$. 

\begin{theorem} \textup{(\cite{14})}
Let $\mathbf V$ be a congruence-modular variety of
algebras, and $A\in\mathbf V$. Let $P$ be a pointed
$A$-overalgebra which is the underlying pointed
$A$-overalgebra of an abelian group overalgebra. Then
there is a unique assignment of abelian group structures
to the pointed sets $_aP$, such that $_a*$ is the zero
element of each $_aP$ and the
functions $\omega^P_{\mathbf a}$ are abelian group
homomorphisms. The group operations in $_aP$ satisfy
$p-p'+p''=d^P_{\langle a,a,a\rangle}(p,p',p'')$.
\end{theorem}

\subsection{Abelian group $A$-overalgebras free on a
pointed $A$-overalgebra}If $\mathbf V$ is
congru\-ence-mod\-u\-lar, a free functor (left
adjoint) for the forgetful functor from $\Ab[A,\mathbf V]$ to $\Pnt[A,\mathbf V]$ is as follows: Given a pointed
$A$-overalgebra $P$, draw the diagram
\[
\CD A @>\iota_P>>
A\semitimes P @>\nat[\kappa,\kappa]>> (A\semitimes
P)/[\kappa,\kappa] \\ @| @VV\pi_PV @VV\pi,V \\ A @= A @= A
\endCD
\]
 where $\kappa=\ker\pi_P$ and $\pi$ is the unique
homomorphism making the diagram commute.
\[
M=\lsil (A\semitimes
P)/[\kappa,\kappa],\pi,\nat[\kappa,\kappa]
\circ\iota_P\rsil
\]
is then a pointed overalgebra, which is an abelian group
overalgebra in a unique way, by theorem~1.1. $M$ is an
abelian group
$A$-overalgebra free on $P$.

 \subsection{$\mathbf V'$
$A$-overalgebras}Suppose $\mathbf V'$ is another variety
of algebras than $\mathbf V$, perhaps of a different type. We
define a $\mathbf V'$ $A$-overalgebra to be an
$A$-overalgebra $M$ such that the $_aM$ are algebras in
$\mathbf V'$ and the $\omega^M_{\mathbf a}$ are homomorphisms
of $\mathbf V'$, and a homomorphism of $\mathbf V'$
$A$-overalgebras to be an $A$-overalgebra homomorphism $f$
such that the $_af$ are homomorphisms of algebras in
$\mathbf V'$. We denote the category of $\mathbf V'$
$A$-overalgebras totally in $\mathbf V$, and homomorphisms
betwen them, by $\mathbf V'[A,\mathbf V]$.

The category $\mathbf V'[A,\mathbf V]$ generalizes
$\Pnt[A,\mathbf V]$ and $\Ab[A,\mathbf V]$ in an obvious way,
and we have $\mathbf V'[A,\mathbf V]=\mathbf V'[\Ov[A,\mathbf V]]$, the category of $\mathbf V'$-objects in the category
$\Ov[A,\mathbf V]$.

Although $\mathbf V'[A,\mathbf V]$ is the category of $\mathbf V'$-objects in the category $\Ov[A,\mathbf V]$, it is also
true that given $M\in\mathbf V'[A,\mathbf V]$, and an $A$-set
$S$, the $A$-functions from $S$ to $M$ form an algebra of
$\mathbf V'$ in a natural way.

If $M\in\mathbf V'[A,\mathbf V]$, and $u$ is an $n$-ary basic
operation, term operation, or element of $\Clo\mathbf V'$,
we will denote by $_au$, or occasionally by $_au^M$, that
operation on the algebra $_aM$.

\subsection{Restriction and induction functors}For each category $\mathbf V'[A,\mathbf V]$, algebra
$X\in\mathbf V$, and homomorphism $f:X\to A$, there is a
functor $_f{\Res:\mathbf V'[A,\mathbf V]\to\mathbf V'[X,\mathbf V]}$, defined by
$
_x({_f{\Res M}})={_{f(x)}M}$
and
$\omega^{{_f{\Res M}}}_{\mathbf x}=\omega^M_{f(\mathbf x)}$,
where $f(\mathbf x)$ stands for $\langle
f(x_1),\ldots,f(x_n)\rangle$.

The restriction functors all have left adjoints, which
are constructed in \cite{14}. We call these \emph{
induction} functors.

We will have occasion to use the functor of induction of
abelian group overalgebras in \S\ref{S:modules}.

\subsection{Products}Let $A$ be an algebra. If $S_1$, $\ldots$, $S_n$ are
$A$-sets, then a product $\Pi_iS_i$ in the category of
$A$-sets is given by
$_a(\Pi_iS_i)={_aS_1}\times\ldots\times{_aS_n}$.
Similarly, if $M_1$, $\ldots$, $M_n\in\mathbf V'[A,\mathbf V]$, then a product $\Pi_iM_i\in\mathbf V'[A,\mathbf V]$ is
given by $_a(\Pi_iM_i)={_aM_1}\times\ldots\times{_aM_n}$,
by 
\[v^{\Pi_iM_i}_{\mathbf a}(\mathbf m_1,\ldots,\mathbf
m_k)=
\langle v^{M_1}_{\mathbf a}(m_{11},\ldots,m_{k1}),\ldots,
v^{M_n}_{\mathbf a}(m_{1n},\ldots,m_{kn})\rangle,
\]
 for
$v\in\Clo_k\mathbf V$, and by
\[
_au(\mathbf
m_1,\ldots,\mathbf m_k)=
\langle{_au^{M_1}(m_{11},\ldots,m_{k1})},\ldots,
{_au^{M_n}(m_{11},\ldots,m_{k1})\rangle},
\]
 for
$u\in\Clo_k\mathbf V'$.

Suppose, on the other hand, that we are given algebras
$A_i\in\mathbf V$, and objects $M_i\in\mathbf V'[A_i,\mathbf V]$,
for $i=1$, $\ldots$, $n$. Let $A=\Pi_iA_i$.
Define the \emph{outer product}
$\boxtimes_iM_i\in\mathbf V'[A,\mathbf V]$ by $_{\mathbf a}(\boxtimes_iM_i)=\Pi_i({_{a_i}M_i})$, and by
\begin{align*}
v^{\boxtimes_iM_i}_{\langle
\mathbf a_1,\ldots,\mathbf a_k\rangle}
(\mathbf
m_1,&\ldots,\mathbf m_k)\\
&=\langle
v^{M_1}_{\langle a_{11},\ldots,
a_{k1}\rangle}(m_{11},\ldots,
m_{k1}),\ldots,v^{M_n}_{\langle a_{1n},\ldots,
a_{kn}\rangle}(m_{1n},\ldots,
m_{kn})\rangle,
\end{align*}
 for $v\in\Clo_k\mathbf V$, and
\[
_{\vec a}u(\mathbf m_1,\ldots,\mathbf m_i)
=\langle{_{a_1}u^{M_1}(m_{11},\ldots,
m_{k1})},\ldots,
{_{a_n}u^{M_n}(m_{1n},\ldots,
m_{kn})}\rangle,
\]
for $u\in\Clo_k\mathbf V'$. That is, $M
=\Pi_i({_{\pi_{A,i}}{\Res M_i}})$, where the
$\pi_{A,i}:A\to A_i$ are the projections to the factors.

\begin{theorem} If all the $A_i$ are the same
algebra $A$, we have
$_{\Delta_A}{\Res(\boxtimes_iM_i)}=M^n$,
where the homomorphism $\Delta_A:A\to\Pi_iA$ is defined by
$a\mapsto\langle a,\ldots,a\rangle$.
\end{theorem}

\subsection{Advantages of the overalgebra
formalism}There are two main advantages for using the formalism of
$A$-overalgebras rather than that of algebras over $A$.
One is that given an object $M\in\mathbf V'[A,\mathbf V]$,
the $_aM$, which are algebras of $\mathbf V'$, are
important objects of study, and a formalism that provides
for this is useful. If we use the formalism of algebras
over $A$, then we will end up considering the $_aM$
anyway, in a different form, as $\pi^{-1}(a)$ for $a\in A$.

The other main advantage is that the formalism of
$A$-overalgebras allows the $_aM$ to be nondisjoint. One
way this helps is to facilitate defining and using the
pointed $A$-overalgebra $\alpha^*$ for $\alpha\in\Con A$.
If we don't allow nondisjoint sets, then we must work
with the algebra often denoted by
$A(\alpha)=A\semitimes\alpha^*$, along with the
projection $\pi_{\alpha^*}$ to $A$.

The other way allowing the $_aM$ to be nondisjoint is
helpful, and this is probably the most important
advantage, is in defining and using the restriction
functors. Using the formalism of algebras over $A$
requires that the restrictions be defined as pullbacks,
and then frequent use must be made of the universal
property of the pullback. This can be done, of course,
but it is tedious, and tends to obscure the situation.

\section{Enveloping Ringoids}\label{S:ringoids}

\subsection{Ringoids}A \emph{ringoid} is a small additive
category. If $\mathbf X$ is a ringoid, with set of objects
$A$, we write
$_{a'}\mathbf X_a$ for $\mathbf X(a,a')$. A left \emph{$\mathbf X$-module} is an additive functor from $\mathbf X$ to $\Ab$,
the category of abelian groups. If $M$ is a left $\mathbf X$-module, we write $_aM$ instead of $M(a)$, and if
$m\in{_aM}$ and $r\in{_{a'}\mathbf X_a}$, we write $rm$
rather than $M(r)(m)$.

If $A$ is an algebra in a variety $\mathbf V$, the \emph{
enveloping ringoid for $A$, with respect to $\mathbf V$}, is
a certain ringoid denoted by $\mathbb Z[A,\mathbf V]$,
which has the underlying set of $A$ as its set of
objects, and such that the category of left $\mathbb
Z[A,\mathbf V]$-modules is isomorphic to the category
$\Ab[A,\mathbf V]$. Enveloping ringoids are treated in
detail in
\cite{14} and
\cite{15}.

\subsection{Construction of the enveloping ringoid}Given an algebra $A$ in a variety $\mathbf V$, we construct
an object $M_a\in\Ab[A,\mathbf V]$ free on an $A$-set
consisting of a singleton at $a$ and the null set
elsewhere, for each $a\in A$. (An abelian group
$A$-overalgebra free on an $A$-set $S$ can be
constructed by constructing an $A$-overalgebra free on
$S$, a pointed $A$-overalgebra free on that, and an
abelian group $A$-overalgebra free on that.)

The enveloping
ringoid can then be given as $_{a'}\mathbb Z[A,\mathbf
V]_a={_{a'}(M_a)}$ for all $a$, $a'\in A$.

\subsection{Abelian group $A$-overalgebras free on an
$A$-set}Once the enveloping ringoid is
defined, there is another method available for defining
abelian group $A$-overlaying algebras free on $A$-sets.
Let $A$ be an algebra in the variety $\mathbf V$, and let
$S$ be an $A$-set. For each $a\in A$, let $_aM$ be the set
of finite formal linear combinations
\[
\sum_ir_is_i,
\]
where $r_i\in{_a\mathbb Z[A,\mathbf V]_{b_i}}$ and
$s_i\in{_{b_i}S}$. This defines an object $M\in\Ab[A,\mathbf V]$ free on $S$.

\section{Abelian
Extensions}\label{S:extensions}

Let $A$ be an algebra of $\mathbf V$, and let
$M\in\Ab[A,\mathbf V]$. We
define an extension in $\mathbf V$ of $A$ by $M$ to be a
triple $\langle\chi,E,\PI\rangle$ such that $E$ is an
algebra of $\mathbf V$, $\PI:E\to A$ is an onto
homomorphism, and $\chi:{_\pi{\Res P}}\to\kappa^*$ is an
isomorphism of pointed overalgebras, where $P$ is the
underlying pointed $A$-overalgebra of $M$ and
$\kappa=\ker\PI$.

For example, given $A\in\mathbf V$, and $M$, totally in
$\mathbf V$, we can form the total algebra $A\semitimes M$.
Let $\kappa=\ker\pi_M$. Let $\chi:{_\pi{\Res
M}}\to\kappa^*$ be the $(A\semitimes M)$-function given by
$_{\pair am}\chi(m')=\pair a{m+m'}$.

\begin{theorem}
$\langle\chi,A\semitimes M,\PI_P\rangle$ is an extension in
$\mathbf V$ of $A$ by $M$. \end{theorem}

\begin{proof} Each $_{\pair am}\chi$ is clearly one-one
and onto, and $\chi$ is an homomorphism of pointed
$(A\semitimes M)$-overalgebras because
$_{\pair am}\chi(0)=\pair am={_{\pair am}*}$, and for
each $v\in\Clo_n\mathbf V$, each $\pair{\mathbf a}{\mathbf m}
=\langle\pair{a_1}{m_1},\ldots,\pair{a_n}{m_n}\rangle\in
(A\semitimes M)^n$, and each $\mathbf m'=\langle
m'_1,\ldots,m'_n\rangle\in{_{\pair{\mathbf a}{\mathbf m}}({_\pi{\Res M}})}$, we have
\begin{align*}
_{v\pair{\mathbf a}{\mathbf m}}\chi
(v^{{_\pi{\Res M}}}_{\pair{\mathbf a}{\mathbf m}}(\mathbf m'))
&={_{v(\pair{\mathbf a}{\mathbf m})}\chi
(v^M_{\mathbf a}(\mathbf m'))}\\
&=\pair{v(\mathbf a)}{v^M_{\mathbf a}(\mathbf m)
+v^M_{\mathbf a}(\mathbf m')}\\
&=v^{A\semitimes M}(\pair{a_1}{m_1+m'_1},\ldots,
\pair{a_n}{m_n+m'_n})\\
&=v^{\kappa^*}_{\pair{\mathbf a}{\mathbf m}}
(\pair{a_1}{m_1+m'_1},\ldots,
\pair{a_n}{m_n+m'_n})\\
&=v^{\kappa^*}_{\pair{\mathbf a}{\mathbf m}}(
{_{\pair{a_1}{m_1}}\chi(m'_1)},\ldots,
{_{\pair{a_n}{m_n}}\chi(m'_n)}
),
\end{align*}
\end{proof}

In this example, the onto homomorphism $\PI=\pi_M$ is split
by the homomorphism $\iota_M$. We say that an
abelian extension $\langle\chi,E,\PI\rangle$ is
\emph{split} if $\PI$ splits.

In general, given an extension $\triple\chi E\pi$ in
$\mathbf V$ of $A$ by $M$, the pointed $E$-overalgebra
$\kappa^*=(\ker\pi)^*$ is isomorphic to the abelian group
overalgebra $_\pi{\Res M}$, and so is itself an abelian
group overalgebra with abelian group operations as given
in theorem~1.1.

\subsection{A lemma using the properties of the
difference term}The following lemma will be very
useful in proving properties of abelian extensions:

\begin{lemma} \label{T:keylemma} If $\triple\chi E\PI$ is
an abelian extension in $\Var$ of $A$ by $M$, and $e$,
$e'$, $e''\in \E$ are such that $\PI(e)=\PI(e')=\PI(e'')$,
then
\[
_e\chi^{-1}(e')+{_{e'}\chi^{-1}(e'')}
={_e\chi^{-1}(e'')}.
\]
\end{lemma}

\begin{proof} Let $\kappa=\ker\PI$, and let $a=\PI(e)$. 
By the properties of $d$ and theorem~1.1, we have
\begin{align*}
{_{e'}\chi^{-1}(e'')}
&=d^M_{\langle
a,a,a\rangle}({_{e'}\chi^{-1}(e'')},{_a0},{_a0})\\
&=d^{_\PI{\Res M}}_{\langle
e',e',e\rangle}
({_{e'}\chi^{-1}(e'')},{_{e'}\chi^{-1}({_{e'}0})},
{_e\chi^{-1}({_e0}))}\\
&={_e\chi^{-1}(d^{\kappa^*}_{\langle e',e',e\rangle}
(e'',{_{e'}0},{_e0}))}\\
&={_e\chi^{-1}(d^E(e'',e',e))}\\
&={_e\chi^{-1}(d^{\makemod\kappa}_{\langle e,e,e\rangle}
(e'',e',{_e0}))}\\
&={_e\chi^{-1}(e'')}-{_e\chi^{-1}(e')}.
\end{align*}
\end{proof}

\subsection{Sections of extensions}Let $\mathcal
E=\triple\chi
\E\PI$ be an abelian extension in $\Var$ of $A$ by
$M$.
Recall that we say that $\mathcal E=\triple\chi
\E\PI$ splits if
there is a homomorphism $\sigma:A\to E$ right inverse to
$\PI$. Not every abelian extension of $A$ by $M$ splits,
as we know, because not every onto homomorphism has a right
inverse.  (We do know that there is always at least one
split extension, given previously.)
At any rate, by the axiom of choice, there is always a
function $\sigma$ such that $\PI\sigma=1_A$, whether or
not the extension $\mathcal E$ splits.

A function $\sigma:A\to E$, not necessarily
a homomorphism, such that $\PI\sigma=1_A$, will be called
a \emph{section} of $\PI$ (or, of the extension
$\mathcal E$).

\subsection{$E$ and $M$}For each $e\in E$,
$_e\chi^{-1}:\PI^{-1}(\PI(e))\to{_{\pi(a)}M}$ is a one-one
and onto function. A section $\sigma$ chooses one
representative $\sigma(a)$ for each $\kappa$-equivalence
class $\PI^{-1}(a)$, and allows us to select functions
$_{\sigma(a)}\chi^{-1}$ which provide an isomorphism of
pointed $A$-sets between $\lsil E,\PI,\sigma\rsil$ and $M$.
We will use this pointed $A$-function extensively in what
follows.

\subsection{Some $A$-sets connected with abelian
extensions}Let $\Var$ be a variety of algebras
of type $\Omega$, and let $A$ be an algebra in $\Var$. We
will define some $A$-sets that will be useful in the study
of abelian extensions of $A$.

Let $X^0_\Var(A)$, or simply $X^0$,
denote the $A$-set
$\lsil A,1_A\rsil$, that is, the underlying set of $A$,
viewed as a set over $A$.  We will refer to the element
$a\in {_aX^0}$, for any $a\in A$, by
$[a]$.  Let $X^1_\Var(A)$, or simply $X^1$,
denote the $A$-set $\lsil S^1,h\rsil$, where $S^1$ is
the set of of pairs (written as follows) $[v;\vec a]$,
where $v\in\Clo_n\Var$ for some
$n$, and $\vec a\in A^n$, and $h[v;\vec a]=v(\vec a)$. 
Let $X^2_\Var(A)$, or simply $X^2$, denote the $A$-set
$\lsil S^2,k\rsil$,
where $S^2$ is the set of triples $[v',\vec v;\vec a]$,
where $\vec a$ is an element of $A^n$ for some $n$, $\vec
v$ is an $n'$-tuple of elements of $\Clo_n\Var$, for
some $n'$, and $v'\in\Clo_{n'}\Var$,
and where $k[v',\vec v;\vec a]=v'(\vec v(\vec a))$.

Note that if $n=0$, $1$, or $2$, then $X^n$ satisfies the
property that if $a\neq a'$, the sets $_aX^n$ and
$_{a'}X^n$ are disjoint. We will often take advantage of
this property in formulas which involve $A$-functions with
domain $X^n$, by dropping the subscript $a$.  Thus, if
$x\in{_aX^n}$ and $\phi$ is an $A$-function from $X^n$ to
another $A$-set, we will write $_a\phi(x)$ simply as
$\phi(x)$.  This is unambiguous because $x$ determines $a$.

\subsection{$\delta^{\mathcal
E}_{\sigma,\sigma'}$}As a first example of an $A$-function with domain one of
the $A$-sets just defined, we define $\delta^{\mathcal
E}_{\sigma,\sigma'}$ to be the $A$-function defined by
\[
\delta^{\mathcal E}_{\sigma,\sigma'}([a])
={_{\sigma(a)}\chi^{-1}(\sigma'(a))},
\]
where $\mathcal E$ is an extension and $\sigma$, $\sigma'$
are two sections for $\mathcal E$.

 \begin{lemma} If $\mathcal E$ is an
abelian extension and $\sigma$ is a
section for $\mathcal E$, then $\delta^{\mathcal
E}_{\sigma,\sigma}\equiv 0$. \end{lemma}

\begin{proof}
$\sigma(a)={_{\sigma(a)}0^{\kappa^*}}$.
Thus, $\delta^{\mathcal
E}_{\sigma,\sigma}={_{\sigma(a)}\chi^{-1}(\sigma(a))}
={_{\sigma(a)}\chi^{-1}(0)}=0$. \end{proof}

\section{Factor Sets}\label{S:factorsets}

\subsection{The factor set of an extension relative to a
section}Let $\mathcal E=\triple\chi E\PI$ be
an abelian extension of $A$ by $M$, and $\sigma$ a
section. The
functions $_{\sigma(a)}\chi^{-1}$, for $a\in A$, allow us
to express the structure of $E$ in terms of $M$.  To begin
we ``measure the failure of $\sigma$ to be a homomorphism''
by defining, for each $[v;\vec a]\in X^1_\Var(A)$, the
element
\begin{align*}
\factset^{\mathcal E,\sigma}[v;\vec a]
&={_{\sigma(v(\vec a))}\chi^{-1}(v(\sigma(\vec
a)))}-{_{\sigma(v(\vec a))}\chi^{-1}(\sigma(v(\vec a)))}\\
&={_{\sigma(v(\vec a))}\chi^{-1}(v(\sigma(\vec
a)))}.
\end{align*}
It is easy to see that $\factset^{\mathcal E,\sigma}$ is
an
$A$-function from $X^1_\Var(A)$ to $M$.
We call $\factset^{\mathcal
E,\sigma}$ the
\emph{factor set for $\mathcal E$, with respect to the
section $\sigma$}.

\begin{theorem} \label{T:splitzero} If $\sigma$ is a
splitting for
$\mathcal E$, then $\factset^{\mathcal E,\sigma}\equiv 0$.
\end{theorem}

\begin{proof} In that case, $\factset^{\mathcal
E,\sigma}[v;\mathbf a]={_{\sigma(v(\vec
a))}\chi^{-1}(v(\sigma(\vec a)))}={_{\sigma(v(\vec
a))}\chi^{-1}(\sigma(v(\vec a)))}={_{v(\vec a)}0}$ for all
$v$ and $\vec a$. \end{proof}

\begin{theorem} Together with $\sigma$, the factor
set $\factset^{\mathcal E,\sigma}$ determines the algebra
structure of $\E$. \end{theorem}

\begin{proof}  Let $v\in\Clo_n\Var$, $\vec e\in E^n$,
and $\vec a=\PI(\vec e)$.  We have
\begin{equation}\label{E:eq1}
\begin{aligned}
_{\sigma(v(\vec
a))}\chi^{-1}(v(\vec e))
&={_{\sigma(v(\vec a))}\chi^{-1}(v_{\sigma(\vec
a)}^{\kappa^*}(\vec e))}\\
&={_{v(\sigma(\vec a))}\chi^{-1}(v_{\sigma(\vec
a)}^{\kappa^*}(\vec e))}
+{_{\sigma(v(\vec a))}\chi^{-1}(v(\sigma(\vec a)))}\\
&=v^{{_\PI{\Res M}}}_{\sigma(\vec a)}({_{\sigma(\vec a)}
\chi
^{-1}(\vec e)}) + \factset^{\mathcal E,\sigma}[v;\vec
a]\\
&=v^M_{\vec a}({_{\sigma(\vec a)}\chi^{-1}
(\vec e)}) + \factset^{\mathcal E,\sigma}[v;\vec a]
.
\end{aligned}
\end{equation}
Thus, if the section
$\sigma$ is given, $v^E$ and $\factset^{\mathcal
E,\sigma}$ determine each other, because the mappings
$_{\sigma(a)}\chi^{-1}$ together make up an isomorphism of
$A$-sets from $\lsil E,\PI\rsil$ to $M$. \end{proof}

\subsection{Abstract factor sets}Because $\E$ is an
algebra in $\Var$, sending each $v\in\Clo_n\Var$ to the
operation $v^\E$, for all $n$, is a clone
homomorphism from $\Clo\Var$ to the clone $\Clo U(E)$ of
all finitary operations on the set $U(E)$.  This means that
for each $\vec v\in(\Clo_n\Var)^{n'}$, and each
$v'\in\Clo_{n'}\Var$, we have ${v'}^E\vec
v^E=(v'\vec v)^E$, where the clone composition on the
left takes place in $\Clo U(E)$, and that on the right
takes place in $\Clo\Var$.  Using equation $(1)$ above,
we obtain, for all $\vec a\in A^n$,
\begin{align*}
\factset^{\mathcal E,\sigma}[v'\mathbf v;\mathbf a]
&={_{\sigma(v'(\vec v(\vec a)))}\chi^{-1}({v'}(\vec
v(\sigma(\vec a))))}\\
&={v'}^M_{\vec
v(\vec a)}({_{\sigma(\vec v(\vec a))}\chi^{-1}(\vec v(\vec
\sigma(\vec a))})+\factset^{\mathcal E,\sigma}[v';\vec
v(\vec a)]\\
&={v'}^M_{\vec v(\vec
a)}\vec v^M_{\vec a}({_{\sigma(\vec a)}\chi^{-1}(\sigma(\vec
a))})+{v'}^M_{\vec v(\vec a)}\factset^{\mathcal
E,\sigma}[\vec v;\vec a]+\factset^{\mathcal
E,\sigma}[v';\vec v(\vec a)]\cr &={v'}^M_{\vec
v(\vec a)}\factset^{\mathcal E,\sigma}[\vec v;\vec
a]+\factset^{\mathcal E,\sigma}[v';\vec v(\vec a)],
\end{align*}
for all
such $\vec v$, $v'$, and $\vec a$, where $\vec v(\vec a)$
stands for $\langle v_1(\vec a),\ldots,v_{n'}(\vec
a)\rangle$. We will call an $A$-function
$\factset:X^1_\Var(A)\to M$ satisfying the family of
equations
\[
\factset[v'\vec v;\vec a]={v'}^M_{\vec v(\vec
a)}\factset[\vec v;\vec a] + \factset[v';\vec v(\vec
a)],
\]
 for each $[v',\vec v;\vec a]\in
X^2_\Var(A)$, a \emph{factor
set for $A$ and $M$}.

Let $\A\in\Var$ and let $M\in\Ab[A,\mathbf V]$.  Given an
arbitrary $A$-function $\factset:X^1_\Var(A)\to M$, we
define an $A$-function $\partial\factset:X^2_\Var(A)\to M$
by the equation
\[
(\partial\factset)[v',\vec v;\vec
a]={v'}^M_{\vec v(\vec a)}\factset[\vec v;\vec
a]-\factset[v'\vec v;\vec a]+\factset[v';\vec v(\vec
a)]
\]
 where $\vec v(\vec a)$ is the $n'$-tuple $\langle
v_1(\vec a),\ldots,v_{n'}(\vec a)\rangle$ and $[\vec v;\vec
a]$ is the $n'$-tuple $\langle [v_1;\vec
a],\ldots,[v_{n'};\vec a]\rangle$.

\begin{theorem} An $A$-function
$\factset:X^1_\Var(A)\to M$ is a factor set for $A$ and
$M$ iff $\partial\factset$ is identically zero, i.e., if
for each $[v',\vec v;\vec a]\in X^2_\Var(A)$, $(\partial
\factset)[v',\vec v;\vec a]={_{{v'}(\vec v(\vec
a))}0}$. \end{theorem}

\subsection{Abstract factor sets and
abelian extensions}
\begin{theorem} \label{T:abstractfact}
Every factor set for $A$ and $M$ arises as
$\factset^{\mathcal E,\sigma}$ for some extension
$\mathcal E$ of
$A$ by $M$ and some section $\sigma$ of that extension.
\end{theorem}

\begin{proof} Given a factor set $\factset$ for $A$
and $M$, we define $E=U(A\semitimes
M)$, and for each $v\in\Clo_n(\Var)$, the
$n$-ary operation $v^E:E^n\to E$, by the equation
\[
v(\pair{a_1}{m_1},\ldots,
\pair{a_n}{m_n})=\pair{v(\vec a)}{v^M_{\vec a}(\vec
m)+\factset[v;\vec a]}.
\]
 The fact that $\factset$ is a
factor set makes the mapping $v\to v^E$ a clone
homomorphism, i.e., makes $E$ an algebra of $\Var$.  For,
let $e_i=\pair{a_i}{m_i}$ for $i=1,\ldots,n$, let $\vec
v$ be an $n'$-tuple of elements of $\Clo_n\mathbf V$, and
let $v'\in\Clo_{n'}\mathbf V$.  We have
\begin{align*}
(v'\vec v)(\vec e)
&=\pair{(v'\vec
v)(\vec a)}{(v'\vec v)^M_{\vec a}(\vec m) +
\factset[v'\vec v;\vec a]}\\
&=\pair{(v'\vec v)(\vec a)}{(v'\vec v)^M_{\vec
a}(\vec m) + {v'}^M_{\vec v(\vec a)}\factset[\vec
v;\vec a]+\factset[v';\vec v(\vec a)]},
\end{align*}
while
\begin{align*}
{v'}(\vec v(\vec e))
&={v'}(\pair{v_1(\vec a)}{(v_1)^M_{\vec a}(\vec
m)+\factset[v_1;\vec a]},\ldots
\pair{v_{n'}(\vec a)}{(v_{n'})^M_{\vec a}(\vec
m)+\factset[v_{n'};\vec a]})\cr
&=\pair{{v'}(\vec v(\vec a))}{{v'}^M_{\vec v(\vec
a)}\vec v^M_{\vec a}(\vec m)+{v'}^M_{\vec v(\vec
a)}\factset[\vec v;\vec a] + \factset[v';\vec v(\vec
a)]},
\end{align*}
 and the desired equality of these two
elements follows from the analogous facts for $A$ and $M$.

 Note $\pi_M:E\to A$ is a
homomorphism with respect to the algebra $E$ just
defined. Let us denote $\ker\pi_M$ by $\kappa$, and
$_{\pi_M}{\Res M}$ by $\bar M$.  We have $\pair
am\mathrel\kappa\pair{a'}{m'}$ iff $a=a'$.  We define
$\chi:\bar M\to\makemod\kappa$ by the equation
\[
_{\pair am}\chi(m')=\pair a{m'+m}\in {_{\pair
am}\kappa^*}.
\]
Each $_{\pair am}\chi$ preserves the
distinguished element because
\[
_{\pair am}\chi({_{\pair am}0})={_{\pair
am}\chi({_a0})}=\pair am={_{\pair am}0}\in{_{\pair
am}\kappa^*}.
\]
Also, $\chi$ preserves the $E$-operations
 because for each
$v\in\Clo_n(\Var)$, $e_i=\pair{a_i}{m_i}\in E$ for
$i=1,\ldots,n$, and $m_i\in
{_{e_i}{\bar M}}$ for $i=1,\ldots,n$, we have
\begin{align*}
v^{\kappa^*}_{\vec e}({_{\vec
e}\chi(\vec{m'})})
&=v^{\kappa^*}_{\vec e}(\pair{a_1}{m'_1+m_1},\ldots,
\pair{a_n}{m'_n+m_n})\\
&=v^E(\pair{a_1}{m'_1+m_1},\ldots,
\pair{a_n}{m'_n+m_n})\\
&=\pair{v(\vec a)}
{v^M_{\vec a}(\vec{m'})+v^M_{\vec a}(\vec
m)+\factset[v;
\vec a]},
\end{align*}
 while on the other
hand,
\begin{align*}
{_{v(\vec e)}\chi(
v^{\bar M}_{\vec e} (\vec{m'}))}
&={_{\langle v(\vec a),v^M_{\vec a}(\vec
m)+\factset[v;\vec a]\rangle}\chi(v^M_{\vec
a}(\vec{m'}))}\\
&=\pair{v(\vec a)}{v^M_{\vec a}(\vec{m'})
+v^M_{\vec a}(\vec
m)+\factset[v;\vec a]}.
\end{align*}
Thus, $\chi$ is a homomorphism of pointed
$E$-overalgebras, and it is clear that each $_e\chi$ is
1-1 and onto.  It follows that $\chi$ is an
isomorphism
of pointed $E$-overalgebras.

Now, let
$\sigma=\iota_\M$ (so that $\sigma(a)=\pair a{{_a0}}$) and
let us show that the factor set we obtain is $\factset$. 
We have
\[
{_{\sigma(\pi_M\pair am)}\chi^{-1}\pair
am}={_{\pair a{{_a0}}}\chi^{-1}\pair am},
\]
 so that the
corresponding factor set is given by
\begin{align*}
\factset^{\mathcal E,\sigma}[v;\vec a]
&=
{_{\sigma(v(\vec a))}\chi^{-1}(v(\sigma(\vec a)))}\\
&=
{_{\sigma(v(\vec
a))}\chi^{-1}(v(\pair{a_1}{{_{a_1}0}},\ldots,
\pair{a_n}{{_{a_n}0}}))}\\ &={_{\pair{v(\vec
a)}{{_{v(\vec a}0}}}\chi^{-1}\pair{v(\vec
a)}{\factset[v;\vec a]}}\\
&=\factset[v;\vec a],
\end{align*}
since by definition,
\[
_{\pair{v(\vec a)}{{_{v(\vec a)}0}}}\chi(\factset[v;\vec
a]) =\pair{v(\vec a)}{\factset[v;\vec a]}.
\]
\end{proof}

\subsection{Effect of choice of section on the corresponding factor
set
}Let us see how factor sets for
the same extension $\mathcal E=\triple\chi{\mathbf E}\PI$ in $\Var$ of $\A$ by $\M$ differ when they are
derived from different sections $\sigma$ and $\sigma'$.
 First, a
definition:

Let $\delta:X^0_\Var(A)\to M$ be an
$A$-function.  We define an
$A$-function $\partial\delta:X^1_\Var(A)\to M$ by the
equation
\[
(\partial\delta)[v;\vec a]=v^M_{\vec a}(\delta[\vec a])
-\delta[v(\vec a)],
\]
for each $[v;\vec a]\in X^1_\Var(\A)$.

\begin{theorem} \label{T:diffboundary} Under these
assumptions, we have
\[
\factset^{\mathcal E,\sigma'}-\factset^{\mathcal
E,\sigma}
=\partial\delta^{\mathcal E}_{\sigma,\sigma'}.
\]
\end{theorem}

\begin{proof} Using equation $(1)$ for $\vec
e=\sigma'(\vec a)$ to expand $\factset^{\mathcal
E,\sigma}[v;\vec a]$, we have
\[
\factset^{\mathcal
E,\sigma'}[v;\vec a]-\factset^{\mathcal E,\sigma}[v;\vec
a] ={_{\sigma'(v(\vec a))}\chi^{-1}(v(\sigma'(\vec a)))}+
v^M_{\vec a}({_{\sigma(\vec a)}\chi^{-1}(\sigma'(\vec
a))})-
{_{\sigma(v(\vec a))}\chi^{-1}(v(\sigma'(\vec a)))}.
\]
But,
\[
{_{\sigma'(v(\vec a))}\chi^{-1}(v(\sigma'(\vec
a)))}={_{\sigma(v(\vec a))}\chi^{-1}(v(\sigma'(\vec a)))}
-{_{\sigma(v(a))}\chi^{-1}(\sigma'(v(\vec a)))}
\]
by lemma 3.1;
combining these two results, we obtain
\begin{align*}
\factset^{\mathcal E,\sigma'}[v;\vec
a]-\factset^{\mathcal E,\sigma}[v;\vec a]
&=v^M_{\vec a}({_{\sigma(\vec a)}\chi^{-1}(\sigma'(\vec
a))})
-{_{\sigma(v(\vec a))}\chi^{-1}(\sigma'(v(\vec
a)))}\\
&=v^M_{\vec
a}(\delta^{\mathcal E}_{\sigma,\sigma'}[\vec a])
-\delta^{\mathcal E}_{\sigma,\sigma'}[v(\vec a)]\\
&=(\partial\delta^{\mathcal E}_{\sigma,\sigma'})[v;\vec
a],
\end{align*}
as was to be proved. \end{proof}

\section{Equivalence of Abelian Extensions}
\label{S:equivalence}

Let $E$ and $\tilde E$ be
algebras of $\Var$, a variety of algebras of type
$\Omega$.  Let $\gamma:E\to\tilde E$ be a
homomorphism, and let $\alpha$ and $\tilde\alpha$ be
congruences of $E$ and
$\tilde E$, respectively, such that
$\gamma(\alpha)\subseteq\tilde\alpha$.  Then we define an
$E$-function,
$\gamma^*:\alpha^*\to{_\gamma{\Res(\tilde\alpha^*)}}$
by the equation
\[_e\gamma^*(e')=\gamma(e').\]

\begin{theorem} $\gamma^*$ is a homomorphism of
pointed $E$-overalgebras. \end{theorem}

\subsection{Equivalent extensions}Let $\mathcal
E=\triple\chi
\E\PI$ and $\tilde{\mathcal E}=\triple{\tilde\chi}{\tilde
\E}{\tilde\PI}$ be abelian extensions in $\Var$ of
$A$ by $M$, where $M\in\Ab[A,\mathbf V]$.  We define an
\emph{equivalence of extensions} from $\mathcal E$ to
$\tilde{\mathcal E}$ to be a homomorphism
$\gamma:E\to\tilde E$, such that
\begin{enumerate}
\item[(1)]
$\PI=\tilde\PI\gamma$, and
\item[(2)]
$\gamma^*\chi={_\gamma{\Res\tilde\chi}}$.
\end{enumerate}

 If $\mathcal E$ and $\hat{\mathcal E}$
are equivalent via an equivalence $\gamma$, we write
$\gamma:\mathcal E\thicksim\hat{\mathcal E}$.

Bearing in mind that,
because of condition (1),
$_\gamma{\Res\circ{_{\tilde\PI}{\Res}}}= {_\PI{\Res}}$, we
express these conditions in the following interrelated
diagrams:
\[
\CD M @. \text{\ \ \ \ } @. {_\PI{\Res M}}
@>\chi>> \makemod\kappa @.\text{\ \ \ \ } @. \E @>\PI>> \A
\\ @| @. @|
@VV\makemod{\gamma}V @. @V\gamma VV @|
\\
M @. \text{\ \ \ \ } @.
{_\PI{\Res M}} @>>{_{\gamma}{\Res\tilde\chi}}>
{_{\gamma}{\Res(\makemod{\tilde\kappa}}}) @.\text{\ \ \ \ }
@. \tilde \E @>>\tilde\PI > \A \\
\\
@. \text{\ \ \ \ } @.
{_{\tilde\PI}{\Res M}} @>>\tilde\chi >
\makemod{\tilde\kappa} \endCD \]
where $\kappa=\ker\pi$ and $\tilde\kappa=\ker\tilde\pi$.

\begin{theorem} Equivalence of extensions is
an equivalence relation on extensions in $\mathbf V$ of $A$
by $M$. \end{theorem}

\begin{proof} If $\mathcal E=\triple\chi E\PI$,
$\tilde{\mathcal E}=\triple{\tilde\chi}{\tilde
E}{\tilde\PI}$, and $\bar{\mathcal
E}=\triple{\bar\chi}{\bar E}{\bar\PI}$ are extensions in
$\mathbf V$ of $A$ by $M$, and $\gamma_1:\mathcal
E\sim\tilde{\mathcal E}$,
$\gamma_2:\tilde{\mathcal E}\sim\bar{\mathcal E}$ are
equivalences, then we first observe that
\[(\gamma_2\gamma_1)^*
=({_{\gamma_1}{\Res\gamma_2^*}})\gamma_1^*.\]
Then, we have
\begin{align*}
(\gamma_2\gamma_1)^*\chi
&=({_{\gamma_1}{\Res\gamma_2^*}})\gamma_1^*\chi\\
&=({_{\gamma_1}{\Res\gamma_2^*}})
{_{\gamma_1}{\Res\tilde\chi}}\\
&={_{\gamma_1}{\Res(\gamma_2^*\tilde\chi)}}\\
&={_{\gamma_1}{\Res({_{\gamma_2}{\Res\bar\chi}})}}\\
&={_{\gamma_2\gamma_1}{\Res\bar\chi}};
\end{align*}
we also have
$\bar\PI\gamma_2\gamma_1=\PI$, whence
$\gamma_2\gamma_1:\mathcal E\sim\bar\mathcal E$.
\end{proof}

\begin{lemma} If $\gamma$ is an equivalence from
$\mathcal E$ to $\tilde{\mathcal E}$, then $\gamma$ is an
isomorphism, and is the unique equivalence from $\mathcal
E$ to $\tilde{\mathcal E}$. \end{lemma}

\begin{proof} $\chi$ and $_\gamma{\Res\chi}$ are
isomorphisms, whence $\gamma^*$ is also by condition (2).
Thus, $\gamma$ maps each $\kappa$-class onto a
$\tilde\kappa$-class, in a one-one fashion. Since each
$\tilde\kappa$-class is the image of a unique
$\kappa$-class by condition (1), these mappings paste
together to give $\gamma$ as an isomorphism.

$\gamma$ is unique, because it is determined
by $\gamma^*$, which is determined by condition
(2). \end{proof}

We denote by
$E_{\Var}(\A,\M)$ the set of equivalence classes
of extensions in $\Var$ of $\A$ by $\M$. It is easy to
see that if $\mathcal E\sim\tilde{\mathcal E}$, then
$\mathcal E$ splits iff $\tilde{\mathcal E}$ splits.
Thus, the split extensions of $A$ by $M$ form one or more
equivalence classes. We will see below, in
corollary~\ref{T:splitequiv}, that all split extensions
are equivalent.

\subsection{Equivalence and factor sets}We wish to relate factor sets and equivalence.

\begin{lemma} \label{T:factequal}
If $\gamma:\mathcal
E\sim\tilde{\mathcal E}$, then given a
section $\sigma$ of $\mathcal E$, we have
$\factset^{\mathcal
E,\sigma}=\factset^{\tilde{\mathcal
E},\gamma\sigma}$. \end{lemma}

\begin{proof}
For each $[v;\vec a]$, we have
\begin{align*}
\factset^{\tilde{\mathcal
E},\gamma\sigma}[v;\vec a] &={_{\gamma(\sigma(v(\vec
a)))}\tilde\chi^{-1}(v(\gamma(\sigma(\vec a))))}\\
&={_{\sigma(v(\vec a))}({_\gamma{\Res\tilde\chi}})^{-1}
({_{\sigma(v(\vec a))}\gamma^*(v(\sigma(\vec a)))})}\\
&={_{\sigma(v(\vec a))}\chi^{-1}(v(\sigma(\vec a)))}\\
&=\factset^{\mathcal E,\sigma}[v;\vec a].
\end{align*}
\end{proof}

\begin{theorem} \label{T:eqdiffbybound} Let $M$ be an
$A$-module.  Let $\mathcal E=\triple\chi\E\PI$ and
$\tilde{\mathcal
E}=\triple{\tilde\chi}{\tilde\E}{\tilde\PI}$ be two
abelian extensions in $\Var$ of $A$ by $M$, with sections
$\sigma$ and $\tau$, respectively.  Then
$\mathcal E$ and $\tilde{\mathcal E}$ are equivalent if
and only if the factor sets $\factset^{\mathcal
E,\sigma}$,
$\factset^{\tilde{\mathcal E},\tau}:X^1_\Var(\A)\to M$
differ by a function of the form $\partial\delta$ where
$\delta:X^0_\Var(A)\to M$ is an $A$-function.
\end{theorem}

\begin{proof} Suppose $\mathcal E$ and
$\tilde{\mathcal E}$ are equivalent via an equivalence
$\gamma$.   By lemma~\ref{T:factequal},
$\factset^{\mathcal E,\sigma}=\factset^{\tilde{\mathcal
E},\gamma\sigma}$.
 On the other hand,
$\factset^{\tilde{\mathcal
E},\gamma\sigma}-\factset^{\tilde{\mathcal
E},\tau}=\partial\delta^{\tilde{\mathcal
E}}_{\tau,\gamma\sigma}$ by theorem~\ref{T:diffboundary}. 
The desired result follows.

Conversely, if the factor set difference
$\factset^{\tilde{\mathcal E},\tau}-\factset^{\mathcal
E,\sigma}=\partial\delta$ for
$\delta:X^0_\Var(A)\to M$ some $A$-function, then we can
produce a new section $\hat\tau$ for $\tilde{\mathcal E}$,
namely
\[\hat\tau(a)={_{\tau(a)}\tilde\chi(\delta[a])},\]
and we have $\factset^{\tilde{\mathcal
E},\hat\tau}=\factset^{\mathcal E,\sigma}$. For,
\begin{align*}
\factset^{\tilde{\mathcal E},\hat\tau}[v;\vec
a] &={_{\hat\tau(v(\vec
a))}\tilde\chi^{-1}(v(\hat\tau(\vec a)))}\\
&={_{\hat\tau(v(\vec a))}\tilde\chi^{-1}(v(
{_{\tau(\vec a)}\tilde\chi(\delta[\vec
a])}))}\\
&={_{\tau(v(\vec a))}\tilde\chi^{-1}(v
({_{\tau(\vec a)}\tilde\chi(\delta[\vec
a])}))}
-{_{\tau(v(\vec a))}\tilde\chi^{-1}(\hat\tau
(v(\vec a)))}\\
&=v^M_{\vec
a}({_{\tau(\vec a)}\chi^{-1}
({_{\tau(\vec a)}\chi(\delta[\vec
a])})})+\factset^{\tilde{\mathcal E},\tau} [v;\vec a]
-{_{\tau(v(\vec a))}\tilde\chi^{-1}
({_{\tau(v(\vec a))}\tilde\chi (\delta[v(\vec
a)]))}}\\
&=\factset^{\tilde{\mathcal
E},\tau}[v;\vec a] +v^M_{\vec a}(\delta[\vec
a])-\delta[v(\vec a)]\\
&=\factset^{\tilde{\mathcal
E},\tau}[v;\vec a] +(\partial\delta)[v;\vec a]\\
&=\factset^{\mathcal E,\sigma}[v;\vec a],
\end{align*}
as desired.

Now we will construct
an isomorphism $\gamma:\E\to\tilde \E$, which is an
equivalence from $\mathcal E$ to $\tilde{\mathcal E}$.  We
define the one-one and onto function $\gamma:E\to\tilde
E$ by
\[\gamma(e)
={_{\hat\tau(\PI(e))}\tilde\chi}({_{\sigma(\PI
(e))}\chi^{-1}(e)}).\]
We have
$\PI=\tilde\PI\gamma$, because $\gamma$ maps each
$\PI^{-1}(a)$ to $\tilde\PI^{-1}(a)$, and we have
$\gamma\sigma=\hat\tau$, because of lemma~3.3.  If we are
given $\vec e\in E^n$, then letting $\vec a=\PI(\vec
e)=\tilde\PI(\gamma(\vec e))$, we have for each
$v\in\Clo_n(\Var)$,
\begin{align*}
\gamma(v(\vec e))
&={_{\hat\tau(v(\vec a))}\tilde\chi(v^M_{\vec
a}({_{\sigma(\vec a)}\chi^{-1}(\vec e)})} +
\factset^{\mathcal E,\sigma}[v;\vec a])\\
&={_{\hat\tau(v(\vec a))}\tilde\chi(v^M_{\vec
a}({_{\hat\tau(\vec a)}\tilde\chi^{-1}(\gamma(\vec e))}) +
\factset^{\tilde{\mathcal E},\hat\tau}[v;\vec a])}\\
&={_{\hat\tau(v(\vec a))}\tilde\chi(
{_{\hat\tau(v(\vec a))}\tilde\chi^{-1}(v(\gamma(\vec
e))})}\\
&=v(\gamma(\vec e)),
\end{align*}
whence the 1-1 and
onto function $\gamma$ is an isomorphism. Finally, for
each $e\in E$, with $a=\PI(e)$, and all $m\in{_aM}$, we
have
\begin{align*}
_e({_\gamma{\Res\tilde\chi}})(m)
&={_{\gamma(e)}\tilde\chi(m)}\\
&={_{\hat\tau(a)}\tilde\chi(
{_{\hat\tau(a)}\tilde\chi^{-1}(
{_{\gamma(e)}\tilde\chi(m)})})}\\
&={_{\hat\tau(a)}\tilde\chi(
{_{\hat\tau(a)}\tilde\chi^{-1}(\gamma(e))}
+{_{\gamma(e)}\tilde\chi^{-1}
({_{\gamma(e)}\tilde\chi(m)})})}\\
&={_{\hat\tau(a)}\tilde\chi({_{\sigma(a)}\chi^{-1}(e)}+m)}
\\
&={_{\hat\tau(a)}\tilde\chi({_{\sigma(a)}\chi^{-1}(e)}
+{_e\chi^{-1}({_e\chi(m)})})}\\
&={_{\hat\tau(a)}\tilde\chi({_{\sigma(a)}\chi^{-1}(
{_e\chi(m)})})}\\
&={_e(\gamma^*\chi)(m)}
;
\end{align*}
thus,
$\gamma^*\chi={_\gamma{\Res\tilde\chi}}$, so
that condition (2) is satisfied, and $\gamma:\mathcal
E\sim\tilde{\mathcal E}$. \end{proof}

\begin{corollary} \label{T:splitequiv} The split
extensions of
$A$ by
$M$ form a single equivalence class. \end{corollary}

\begin{proof} If $\mathcal E$ is an extension in $\mathbf V$
of
$A$ by $M$, and $\sigma$ is a splitting, then
$\factset^{\mathcal E,\sigma}\equiv 0$ by
theorem~\ref{T:splitzero}. If
$\tau$ is another section of $\mathcal E$, then
$\factset^{\mathcal E,\tau}=\partial\delta$ for some
$A$-function $\delta:X^0\to M$, by
theorem~\ref{T:diffboundary}. Thus, every factor set of a
split extension has the form
$\partial\delta$. It follows by the theorem that all
split extensions are equivalent. \end{proof}

\section{$\mathbf E_{\mathbf V}(A,M)$ as a
Homology Object}\label{S:homologyobject}

\subsection{A fragment of a cochain
complex}We previously introduced $A$-sets
$X^0=X^0_{\mathbf V}(A)$, $X^1=X^1_{\mathbf V}(A)$, and
$X^2=X^2_{\mathbf V}(A)$, and, with the additional data of an
object $M\in\Ab[A,\mathbf V]$, abelian group homomorphisms
$\partial$, so that the diagram
\[A\hbox{-}\Set(X^0,M)\overset\partial\to A\hbox{-}
\Set(X^1,M)\overset\partial\to A\hbox{-} \Set(X^2,M)\]
is a fragment of a cochain complex in the abelian
category $\Ab$.

\subsection{The cohomology group}An $A$-function $\factset:X^1\to M$ is a factor set for
$A$ and $M$ iff $\partial\factset = 0$, by definition. We
showed (theorem~\ref{T:abstractfact}) that each such
abstract factor set is a factor set $\factset^{\mathcal
E,\sigma}$ for some extension
$\mathcal E$ and section $\sigma$. On the other hand, two
factor sets $\factset$, $\factset'$ correspond to
equivalent extensions iff
$\factset-\factset'=\partial\delta$ for some $A$-function
$\delta:X^0\to M$. It follows that

\begin{theorem} The
equivalence classes of extensions of $A$ by $M$, i.e.,
the elements of the set $\mathbf E_{\mathbf V}(A,M)$,
correspond to the elements of the cohomology group of the
above fragment of a cochain complex in the abelian
category $\Ab$.
\end{theorem}

\subsection{The equivalence class of split
extensions}
\begin{theorem} The class
of split extensions is the zero element of $\mathbf E_{\mathbf V}(A,M)$.
\end{theorem}

\begin{proof} Choosing a split extension $\mathcal E$,
and a splitting $\sigma$ as the section, we obtain
$\factset^{\mathcal E,\sigma}\equiv 0$, by
theorem~\ref{T:splitzero}.
\end{proof}

 \subsection{$\mathbf E_{\mathbf V}(A,M)$ as an
abelian group}The correspondence of
theorem~6.1 allows us to place the structure of an
abelian group on the set $\mathbf E_{\mathbf V}(A,M)$:
If $\mathcal E_1$, $\mathcal E_2$ are represented by
factor sets $\factset_1$ and $\factset_2$, then $[\mathcal
E_1]+[\mathcal E_2]$ is represented by
$\factset_1+\factset_2$.

\begin{theorem} If $\mathcal
E=\triple\chi E\PI$ is an extension of $A$ by $M$, then
$-[\mathcal E]=[\triple{-\chi}E\PI]$.
\end{theorem}

\begin{proof} $-\chi$ is also an isomorphism, and
the formula for factor sets shows that
\[\factset^{\triple{-\chi}E\pi,\sigma}
=-\factset^{\triple\chi
E\pi,\sigma},\]
for any section $\sigma$.
\end{proof}

\section{Compositions of Extensions and
Homomorphisms}\label{S:composition}

For this section, let $\mathcal E=\triple\chi E\PI$ be an
abelian extension of $A$ by $M$. We will describe ways of
composing $\mathcal E$ with homomorphisms in $\mathbf V$ and
$\Ab[A,\mathbf V]$. Conceptually, these composition
operations are best seen as operations on the \emph{
extension class} $[\mathcal E]$, the equivalence class of
extensions in $\mathbf V$ of $A$ by $M$ represented by
$\mathcal E$. However, it is not hard to define the
composition of an extension and an algebra homomorphism as
a specific extension, and sometimes useful, so we will
give the definition of that composition in that form.

\subsection{The extension $\mathcal Eg$ for a homomorphism
$g:A'\to A$}Suppose that we have
another algebra $A'$ and a homomorphism $g:A'\to
A$. We will show how to create from $\mathcal E$ an
extension of $A'$ by $_g{\Res M}$. The first step is to
consider the fibered product $E'$ of $A'$ and $E$, with
the associated homomorphisms
which we label $\PI'$ and $\gamma$:
\[\CD E' @>\PI'>> A' \\
@V\gamma VV @VVg.V \\
E @>>\PI>A\endCD \]

\begin{lemma} In this situation, if
$\kappa=\ker\PI$ and $\kappa'=\ker\PI'$, then
$\gamma^*:{\kappa'}^*\to{_\gamma{\Res\kappa^*}}$ is an
isomorphism of pointed $E'$-overalgebras. \end{lemma}

\begin{proof}
If $\pair e{a'}$ is an element of $E'$, then the elements
of $_{\pair e{a'}}{\kappa'}^*$ are the elements of $E'$
of the form
$\pair {e'}{a'}$, and such a pair will belong to $E'$ iff
$\PI(e')=g(a')$. On the
other hand, $_{\pair
e{a'}}({_\gamma{\Res\kappa^*}})={_e\kappa^*}$ is the set
of $e'$ such that $\PI(e)=\PI(e')$.
However, we have $g(a')=\PI(e)$, because $\pair e{a'}\in
E'$.  Thus, the mapping $\pair{e'}{a'}\mapsto e'$ is
one-one and onto. Since this applies to each $_{\pair
e{a'}}\gamma^*$, $\gamma^*$ is an
isomorphism. \end{proof}

Now, we define
$\chi'=(\gamma^*)^{-1}\circ{_\gamma{\Res\chi}}$. Since
$\gamma^*$ and $\chi$ are isomorphisms, so is $\chi'$.
Also, we have $\PI\gamma=g\PI'$, so that
$\chi':{_{\PI'}{\Res({_g{\Res M}})}}\cong{\kappa'}^*$.
Thus, we obtain an extension $\mathcal
Eg=\triple{\chi'}{E'}{\PI'}$ of $A'$ by $_g{\Res M}$. We
will denote
the extension class $[\mathcal Eg]$ by $[\mathcal E]g$.

If $\sigma$ is any section for $\mathcal E$, let
$\sigma':A'\to E'$ be defined by
$\sigma'(a')=\pair{\sigma(g(a'))}{a'}$. Then, we have

\begin{lemma} If $\factset=\factset^{\mathcal
E,\sigma}$, then there is an extension $\bar{\mathcal
E}\in[\mathcal E]g$, and a section $\bar\sigma$, such that
or each $[v;\vec c\,]\in X^1_{\mathbf V}(A')$,
$\factset^{\bar{\mathcal E},\bar\sigma}[v;\vec
c\,]=\factset[v;g(\vec c\,)]$. \end{lemma}

\begin{proof} Let $\bar{\mathcal E}=\mathcal Eg$ and
$\bar\sigma=\sigma'$ as defined above.
We have
\begin{align*}
\factset^{\mathcal Eg,\sigma'}[v;\vec c\,]
&={_{\sigma'(v(\vec c))}{\chi'}^{-1}(v(\sigma'(\vec
c\,)))}\\
&={_{\sigma'(v(\vec
c))}({\gamma^*}^{-1}\circ{_\gamma{\Res
\chi}})^{-1}(v(\sigma'(\vec c\,))}\\
&={_{\sigma'(v(\vec
c))}({_\gamma{\Res\chi^{-1}}}) (\gamma^*(v(\sigma'(\vec
c\,))))}\\
&={_{\sigma(v(g(\vec c)))}\chi^{-1}(v(\sigma(g(\vec
c\,))))}\\
&=\factset^{\mathcal E,\sigma}[v;g(\vec c\,)].
\end{align*}
\end{proof}

\begin{corollary} \label{T:Eghassoc} If in addition to
$\mathcal E$,
$A'$, and $g$ we have an algebra $A''$ and homomorphism
$h:A''\to A'$, then $\mathcal E(gh)\sim (\mathcal Eg)h$.
\end{corollary}

\begin{corollary} \label{T:Eghom}
The mapping $[\mathcal E]\mapsto[\mathcal E]g$ is
a group homomorphism from $\mathbf E_{\mathbf V}(A,M)$ to
the group
$\mathbf E_{\mathbf V}(A',{_g{\Res M}})$.
\end{corollary}

\subsection{The extension class $\dot g[\mathcal E]$ for a
homomorphism $\dot g:M\to M'$}Suppose again
that $\mathcal E=\triple\chi E\PI$ is given, and that $M'$
is another object of $\Ab[A,\mathbf V]$, and $\dot g:M\to
M'$ a homomorphism. $\mathcal E$ determines an extension
class
$[\mathcal E]$. We will specify an
extension class $\dot g[\mathcal E]$ of $A$ by $M'$ by
giving a representative factor set.

Let $\factset^{\mathcal E,\sigma}$ be a factor set for
$\mathcal E$, where $\sigma$ is some section. Since
$\factset^{\mathcal E,\sigma}$ is an $A$-function from
$X^1$ to $M$, and $\dot g$ is a homomorphism from $M$ to
$M'$, we can compose them and obtain an $A$-function
$\factset'=\dot g\factset^{\mathcal E,\sigma}:X^1\to M'$.

\begin{theorem} $\factset'$ is a factor set
for $A$ and $M'$.\end{theorem}

\begin{proof}
We have
\begin{align*}
(\partial\factset')[v',\vec v;\vec a]
&={v'}^{M'}_{\vec v(\vec a)}\factset'[\vec v;\vec a]
-\factset'[v'\vec v;\vec a]+\factset'[v';\vec v(\vec
a)]\\
&=\dot g{v'}^M_{\vec v(\vec a)}\factset[\vec v;\vec a]-
\dot g\factset[v'\vec v;\vec a]
+\dot g\factset[v';\vec v(\vec a)]\\
&=\dot g(\partial\factset)[v',\vec v;\vec a]\\
&=\dot g({_{v'\vec v(\vec a)}0^M})\\
&={_{v'\vec v(\vec a)}0^{M'}}.
\end{align*}
 \end{proof}

 Thus, by the construction of
theorem~\ref{T:abstractfact},
$\factset'$ is the factor set of some extension, which we
shall for the moment denote by $\mathcal E^{\sigma,\dot
g}$.

\begin{theorem} If $\sigma'$ is another section
of $\mathcal E$, then $\mathcal E^{\sigma,\dot
g}\sim\mathcal E^{\sigma',\dot g}$. I.e., $[\mathcal
E^{\sigma,\dot g}]$ does not depend on $\sigma$.
\end{theorem}

\begin{proof} $\factset^{\mathcal
E,\sigma}-\factset^{\mathcal E,\sigma'}=\partial\delta$
for some $A$-function $\delta:X^0\to M$. Then,
$\dot g\factset^{\mathcal E,\sigma}-\dot
g\factset^{\mathcal E,\sigma'}=\dot g(\partial\delta)$.
However,
$\dot g(\partial\delta)=\partial(\dot g\delta)$. For,
\begin{align*}
\dot g(\partial\delta)[v;\vec a]
&=\dot g(v^M_{\vec a}(\delta[\vec a])-\delta[v(\vec
a)])\\
&=v^{M'}_{\vec a}(\dot g\delta[\vec a])-\dot
g\delta[v(\vec a)]\\
&=\partial(\dot g\delta)[v;\vec a];
\end{align*}
Thus, $\mathcal E^{\sigma,\dot g}\sim\mathcal
E^{\sigma',\dot g}$ by theorem~\ref{T:eqdiffbybound}.
\end{proof}

We denote $[\mathcal E^{\sigma,\dot g}]$ by $\dot
g[\mathcal E]$.

\begin{theorem} \label{T:hgEassoc}
If $M''$ is another
object of
$\Ab[A,\mathbf V]$, and $\dot h:M'\to M''$ is another
homomorphism, then $(\dot h\dot g)[\mathcal E]=\dot
h(\dot g[\mathcal E])$. \end{theorem}

\begin{theorem} \label{T:gEhom}
The mapping $[\mathcal
E]\mapsto\dot g[\mathcal E]$ is an abelian group
homomorphism from $\mathbf E_{\mathbf V}(A,M)$ to $\mathbf E_{\mathbf V}(A,M')$.
\end{theorem}

\begin{proof} Follows from the definition of $\dot g$
being a homomorphism of abelian group objects in the
category $\Ov[A,\mathbf V]$. \end{proof}

\subsection{Relationship of these two compositions}
\begin{theorem} \label{T:gEgassoc} If
$\mathcal E$ is an extension of $A$ by $M$, and $g$ and
$\dot g$ are given, then
$(\dot g[\mathcal E])g=({_g{\Res\dot g}})([\mathcal E]g)$.
\end{theorem}

\begin{proof} We will show that representative factor sets
for these two extension classes are equal.

Let $\factset$ be a factor set representing the homology
class corresponding to the extension class of $\mathcal
E$. Then $\dot g\factset$ represents $\dot g[\mathcal E]$.
A factor set $\factset'$ for $(\dot g[\mathcal E])g$ can
then be defined by
\[\factset':[v;\vec c\,]\mapsto(\dot g\factset)[v;g(\vec
c\,)].\]
On the other hand, a factor set $\bar\factset$ for
$[\mathcal E]g$ can be defined by
\[\bar\factset:[v;\vec c\,]\mapsto\factset[v;g(\vec
c\,)]\] and then a factor set for $({_g{\Res\dot
g}})([\mathcal E]g)$ can be given by $({_f{\Res\dot
g}})\bar\factset$. For each $[v;\vec c\,]$, we have
\begin{align*}
({_g{\Res\dot g}})\bar\factset[v;\vec c\,]
&={_{g(v(\vec c\,))}\dot g(\factset[v;g(\vec c\,)])}\\
&=\dot g\factset[v;g(\vec c\,)]\\
&=\factset'[v;\vec c\,].
\end{align*}
\end{proof}

\section{The Addition Operation
on $\mathbf E_{\mathbf V}(Q,M)$}\label{S:addition}

Now we will describe a method for constructing
the result of adding a pair $\mathcal E_1$, $\mathcal E_2$
of extensions of $A$ by $M$. This method closely parallels
the method of adding two module extensions called the Baer
sum. However, while the Baer sum uses only universal
properties to construct an extension having the desired
properties, we can do this only up to a point, and must
complete the construction of a representative of
$[\mathcal E_1]+[\mathcal E_2]$ by using factor sets and
theorem~\ref{T:abstractfact}.

\subsection{Outer product of a finite
tuple of extensions}Let $\mathcal
E_i=\triple{\chi_i}{E_i}{\PI_i}$ be an extension of $A_i$
by $M_i$ for $i=1$, $\ldots$, $n$, where $A_i\in\mathbf V$
and $M_i\in\Ab[A_i\mathbf V]$ for all $i$. Define
$E=\Pi_iE_i$ and $A=\Pi_iA_i$. Define
$\PI:E\to A$ by the equation $\PI(\mathbf e)=\langle\PI_1(e_1),\ldots,\PI_n(e_n)\rangle$. We will
construct an extension $\mathcal E=\triple\chi E A$ of $A$
by
$\boxtimes_iM_i$ from the extensions $\mathcal E_i$.

If $\kappa_i=\ker\PI_i$ for all $i$,
$\kappa=\ker\PI$, and 
$\pi_{E,i}:E\to E_i$ are the projections to the
factors, then we have

\begin{lemma} $\kappa^*$ is naturally
isomorphic to $\Pi_i({_{\pi_{E,i}}{\Res\kappa_i^*}})$.
\end{lemma}

We define $\chi:{_\PI{\Res
M}\to\kappa^*}$ by the equation \[_{\mathbf e}\chi(\mathbf
m)=\langle{_{e_1}\chi_1(m_1)},\ldots,
{_{e_n}\chi_n(m_n)}\rangle,\] and finally, we have

\begin{theorem} $\triple\chi E\PI$ is an
extension of $A$ by $M$. \end{theorem}

We denote this extension by $\boxtimes_i\mathcal E_i$.

\subsection{A factor set of $\boxtimes_i\mathcal E_i$}Let
$\sigma_i$ be a section of $\mathcal E_i$ for each $i$,
and let $\sigma:A\to E$ be the section of
$\boxtimes_i\mathcal E_i$ defined by $\sigma(\mathbf a)=\langle\sigma_1(a_1),\ldots,\sigma_n(a_n)\rangle$.

\begin{theorem} For each $[v;\langle\mathbf a_1,\ldots,\mathbf a_\ell\rangle]\in X^1_{\mathbf V}(A)$, we
have
\[
\factset^{\boxtimes_i\mathcal
E_i,\sigma}[v;\langle\mathbf a_1,\ldots,\mathbf
a_\ell\rangle]=\langle
\factset^{\mathcal E_1,\sigma_1}[v;\langle a_{11},\ldots,
a_{\ell 1}\rangle],\ldots,\factset^{\mathcal E_n,\sigma_n}
[v;\langle a_{1 n},\ldots,a_{\ell n}\rangle]\rangle.
\]
\end{theorem}

\begin{proof} We have
\begin{align*}
\factset^{\boxtimes_i\mathcal
E_i,\sigma}&[v;\langle \mathbf a_1,\ldots,\mathbf a_\ell\rangle]\cr
&={_{\sigma(v(\mathbf a_1,\ldots,\mathbf
a_\ell))}
\chi^{-1}(v(\sigma(\mathbf
a_1),\ldots,\sigma(\mathbf a_\ell)))}\\
&=\langle{_{\sigma_1(v(a_{11},\ldots,a_{\ell
1}))}\chi_1^{-1}(v(\sigma_1(a_{11},\ldots,a_{\ell
1})))},\\
&\qquad\qquad\ldots,
{_{\sigma_n(v(a_{1n},\ldots,a_{\ell
n}))}\chi_1^{-1}(v(\sigma_n(a_{1n},\ldots,a_{\ell
n})))}\rangle\\
&=\langle\factset^{\mathcal
E_1,\sigma_1}[v;\langle a_{11},\ldots,a_{\ell
1}\rangle],\ldots,
\factset^{\mathcal E_n,\sigma_n}[v;\langle a_{1n},\ldots,
a_{\ell n}\rangle]\rangle.
\end{align*}
\end{proof}

\subsection{The addition operation on $\mathbf E_{\mathbf V}(Q,M)$}We observe that $+^M:M^2\to M$ is an arrow of
$\Ab[A,\mathbf V]$. Given a pair $\pair{\mathcal
E_1}{\mathcal E_2}$ of extensions in $\mathbf V$ of $A$ by
$M$,
$(\mathcal E_1\boxtimes\mathcal E_2)\Delta_A$ is an
extension of
$A$ by $_{\Delta_A}{\Res M^{\boxtimes 2}}=M^2$. Thus, we
can form the composite $+^M[(\mathcal
E_1\boxtimes\mathcal E_2)\Delta_A]$, which is a class of
extensions in $\mathbf V$ of $A$ by $M$.

\begin{theorem} We have
\[+^M[(\mathcal E_1\boxtimes\mathcal
E_2)\Delta_A]=[\mathcal E_1]+[\mathcal E_2].\]
\end{theorem}

\begin{proof}
Let sections $\sigma_1$, $\sigma_2$ of the extensions
$\mathcal E_i$ be chosen. The factor set
$\factset^{\mathcal E_1,\sigma_1}+\factset^{\mathcal
E_2,\sigma_2}$ is a factor set of $[\mathcal
E_1]+[\mathcal E_2]$. On the other hand, a factor set for
$+^M[(\mathcal E_1\boxtimes\mathcal E_2)\Delta_A]$ can be
computed from the factor set
$\factset^{\mathcal E_1\boxtimes\mathcal E_2,\sigma}$
computed in theorem 8.3. It is easy to see that these
factor sets are equal. \end{proof}

\section{Abelian Extensions as a Bifunctor}
\label{S:bifunctor}

We want to be able to treat the $\mathbf W$-algebra of
abelian extensions of $A$ by $M$ as a bifunctor. On
objects, we define the bifunctor
\[\mathbf E_{\mathbf V}:(\Ov[A,\mathbf
V])^{\text{op}}\times\Ab[A,\mathbf V]\to\Ab\] by the
formula
\[\mathbf E_{\mathbf V}(Q,M)=\mathbf E_{\mathbf V}
(A\semitimes Q,{_{\PI_Q}{\Res M}}),\]
where we recall that $\pi_Q:A\semitimes Q\to A$ is defined
by $\pi_Q:\pair aq\mapsto a$.

On arrows, it suffices to define an abelian
group homomorphism $\mathbf E_{\mathbf V}(r,1_M)$,
which we will write as a composition $[\mathcal
E]\mapsto[\mathcal E]r$, and an abelian group homomorphism
$\mathbf E_{\mathbf V}(1_Q,\dot g)$, which we will
write as a composition on the other side, $[\mathcal
E]\mapsto\dot g[\mathcal E]$, and to prove the following:
\begin{enumerate}
\item $([\mathcal E]r)s=[\mathcal E](rs)$;
\item $\dot h(\dot g([\mathcal E]))=(\dot h\dot
g)[\mathcal E]$; and
\item $\dot g([\mathcal E]g)=(\dot g[\mathcal E])g$,
\end{enumerate}
when those compositions are defined.

 If
$\mathcal E\in\mathbf E_{\mathbf V}(Q,M)$,
$Q'\in\Ov[A,\mathbf V]$, and $r:Q'\to Q$ is a homomorphism,
then we define $[\mathcal E]r=[\mathcal E](A\semitimes
r)$. 
$[\mathcal E]r$ is an element of
\begin{align*}
\mathbf E_{\mathbf V}(A\semitimes
Q',{_{A\semitimes r}{\Res({_{\pi_Q}{\Res M}})}})
&=\mathbf E_{\mathbf V}(A\semitimes
Q',{_{\pi_{Q'}}{\Res M}})\\
&=\mathbf E_{\mathbf V}(Q',M).
\end{align*}
From corollary~\ref{T:Eghom}, this mapping is an abelian
group homomorphism, and it is easy to see that property
(1) holds, from the corresponding fact
(corollary~\ref{T:Eghassoc}) for composition of extensions
with homomorphisms of
$\mathbf V$.

On the other hand, if $M'$ is another object of
$\Ab[A,\mathbf V]$ and $\dot g:M\to M'$ is a homomorphism, we
define \[\dot g[\mathcal E]=({_{\pi_Q}{\Res \dot
g}})[\mathcal E].\] This is an element of $\mathbf
E_{\mathbf V}(Q,M')$ and it is clear from
theorem~\ref{T:gEhom} that the mapping is an abelian group
homomorphism, and from theorem~\ref{T:hgEassoc} that
property (2) holds.

Finally, if $\mathcal E$, $Q'$, $r$, $M'$, and $\dot g$
are all given, we must prove that $(\dot g[\mathcal
E])r=\dot g([\mathcal E]r)$, and this will complete the
proof that
$\mathbf E_{\mathbf V}$ is a bifunctor as we have
defined it. Using theorem~\ref{T:gEgassoc}, we have
\begin{align*}
(\dot g[\mathcal E])r
&=(({_{\pi_Q}{\Res\dot g}})[\mathcal E])(A\semitimes r)\\
&={_{A\semitimes r}{\Res({_{\pi_Q}{\Res\dot
g}})}}([\mathcal E](A\semitimes r))\\
&={_{\pi_{Q'}}{\Res\dot g}}([\mathcal E]r)\\
&=\dot g([\mathcal E]r).
\end{align*}

\section{Module Extensions}\label{S:modules}

\subsection{$R$-modules as a variety of algebras}If $R$ is a ring, then the left $R$-modules can be
treated as algebras having an underlying abelian group
structure and one unary operation for each element of $R$.
The reader can easily supply a list of identities
defining the variety of left $R$-modules. This is a
congruence-modular variety, as can easily be proved,
because of the underlying abelian group structures.

It is important to note for what follows that if
$v\in\Clo_n R$-$\Mod$, then $v(\mathbf m)$ has the form
$\Sigma_ir_im_i$.

\subsection{Abelian group overalgebras in the variety of
left $R$-modules}We need to prove some facts about abelian group
overalgebras totally
in $R$-$\Mod$:

 \begin{theorem} Let $A$ and $B$ be
left $R$-modules, and $f:A\to B$ an onto homomorphism. Then
every abelian group $A$-overalgebra totally in $R$-$\Mod$
is isomorphic to the restriction, along $f$, of an
abelian group $B$-overalgebra.
\end{theorem}

\begin{proof} Let $M$ be an abelian group $A$-overalgebra
totally in $R$-$\Mod$, and let $\pair{M'}\eta$ be a
universal arrow to the functor $_f{\Res}$. That is, $M'$ is
an induced abelian group $B$-overalgebra in $R$-$\Mod$ of
$M$ along $f$.
By \cite{14, theorem~10.6}, since $f$ is onto, we can
equally well consider $M'$ to be an induced pointed
$B$-overalgebra of $M$ along $f$. Thus,
\cite{14, theorem~C.6.5} applies, and the diagram
\[\CD A\semitimes M @<<\iota_M< A \\
@V\mu(A\semitimes\eta_M)VV @VVfV \\
B\semitimes M' @<<\iota_{M'}< B \endCD \]
is a pushout diagram, where $\mu$ is the \emph{mashing
homomorphism } defined by $\mu:\pair
ax\mapsto\pair{f(a)}x$. However, $\pi_M\iota_M=1_A$,
whence $A\semitimes M\cong A\oplus{_0M}$. Similarly,
$B\semitimes M'\cong B\oplus{_0M'}$. It follows easily that
$_0M\cong{_0M'}$, and that $M\cong{_f{\Res M'}}$.
\end{proof}

\begin{corollary} If $M$ is an abelian group
$A$-overalgebra totally in $R$-$\Mod$, then $M$ is
isomorphic to a restriction, along the unique
homomorphism $\pi_0:A\to 0$, of an abelian group
$0$-overalgebra. \end{corollary}

\begin{corollary} \label{T:determined} If $M$, $M'$ are
abelian group
$A$-overalgebras totally in $R$-$\Mod$, and $\chi:M\to M'$
is a homomorphism, then $\chi$ is determined by $_0\chi$.
\end{corollary}

\begin{proof} By \cite{14, theorem~C.6.4} restriction
along $\pi$ is a full functor from $\Pnt[0,R$-$\Mod]$ to
$\Pnt[A,R$-$\Mod]$. Since $M$ and $M'$ are isomorphic to
restrictions along $\pi$,
$\chi=\phi^{-1}\circ({_\pi{\Res\hat\chi}})\circ\phi'$,
where $\phi:M\cong{_\pi{\Res \hat M}}$,
$\phi':M'\cong{_\pi{\Res\hat M'}}$, and $\hat\chi:\hat
M\to\hat M'$. It follows that $\chi$ is determined by any
of its components $_a\chi$. \end{proof}

\subsection{Abelian extensions and module
extensions}
\begin{theorem} Let $Q\in\Ov[0,R$-$\Mod]$, and let
$M\in\Ab[0,R$-$\Mod]$. Then $\mathbf E_{\mathbf V}(Q,M)\cong\Ext(0\semitimes Q,{_0M})$.
\end{theorem}

\begin{proof}
Given an extension $\triple\chi E\pi$ of $0\semitimes Q$
by $_{\pi_Q}{\Res M}$, let $\iota={_0\chi}$, and view
$\iota$ as a homomorphism from $_0M$ to $E$. Then we have
a module extension
\[0\to{_0M}\overset\iota\to E\overset\pi\to
0\semitimes Q\to 0\]
of $0\semitimes Q$ by $_0M$.

On the other hand, given $\iota$, $E$, and $\pi$, we
define $_e\chi(m)=e+\iota(m)$. This gives an $E$-function
which is clearly one-one and onto, and is an
$E$-overalgebra homomorphism from $_{\pi_Q\pi}{\Res M}$
to $\kappa^*$, where $\kappa=\ker\pi$, because
\begin{align*}
v_{\mathbf e}^{\kappa^*}({_{\mathbf e}\chi(\mathbf m)})
&=v^{\kappa^*}_{\mathbf e}(e_1+\iota(m_1),
\ldots,e_n+\iota(m_n))\\
&=v(e_1+\iota(m_1),\ldots,e_n+\iota(m_n))\\
&=v(\mathbf e)+v(\iota(\mathbf m))\\
&=v(\mathbf e)+\iota(v^M(\mathbf m))\\
&=v(\mathbf e)+\iota(v^M_{\langle 0,\ldots,0\rangle}
(\mathbf m))\\
&={_{v(\mathbf e)}\chi(v^M_{\langle
0,\ldots,0\rangle}(\mathbf m))}\\
&={_{v(\mathbf e)}\chi(v^{{_{\pi_Q\pi}
{\Res M}}}_{\mathbf e}(\mathbf
m))};
\end{align*}
thus, $\chi$ is an isomorphism, and
$\triple\chi E\pi$ is an extension in $R$-$\Mod$ of
$0\semitimes Q$ by
$_{\pi_Q}{\Res M}$, i.e., an element of $\mathbf E_{\mathbf V}(Q,M)$.

The first construction, of $\iota$ from $\chi$, is
one-one by corollary~\ref{T:determined}, but it is also
onto, because if we start with $\iota$ and construct
$\chi$, then we get $\iota$ back as $\iota={_0\chi}$.
Thus, the two mappings are inverses to each other.

If two extensions $\mathcal E=\triple\chi E\pi$ and
$\tilde{\mathcal E}=\triple{\tilde\chi}{\tilde
E}{\tilde\pi}$ in $R$-$\Mod$ of $0\semitimes Q$ by $M$
are given, and
$\gamma:\mathcal E\sim\tilde{\mathcal E}$, then from
$\gamma^*\chi={_\gamma{\Res\tilde\chi}}$ we obtain that
$\gamma{_0\chi}={_0\gamma^*}{_0\chi}={_0\tilde\chi}$,
whence we have a commutative diagram
\[\CD {_0M} @>_0\chi>> E @>>\pi> 0\semitimes Q \\
@| @V\gamma VV @| \\
{_0M} @>>{_0\tilde\chi}> \tilde E @>>\tilde\pi>
0\semitimes Q\endCD\]
showing that the constructed module extensions remain
equivalent, by the customary definition. On the other
hand, given a commutative diagram
\[\CD {_0M} @>\iota>> E @>>\pi> 0\semitimes Q \\
@| @V\gamma VV @| \\
{_0M} @>>\tilde\iota> \tilde E @>>\tilde\pi>
0\semitimes Q\endCD
\]
then let $\chi$ and $\tilde\chi$ be constructed from
$\iota$ and $\tilde\iota$: we have $\tilde\pi\gamma=\pi$
and for all $e\in E$, and $m\in{_0M}$,
\begin{align*}
_e\gamma^*{_e\chi(m)}
&={_e\gamma^*(e+\iota(m))}\\
&=\gamma(e)+\gamma\iota(m)\\
&=\gamma(e)+\tilde\iota(m)\\
&={_{\gamma(e)}\tilde\chi(m)}\\
&={_e({_\gamma{\Res\tilde\chi}})(m)},
\end{align*}
whence $\gamma^*\chi={_\gamma{\Res\tilde\chi}}$. Thus,
the constructed elements of $\mathbf E_{\mathbf V}(Q,M)$ are
equivalent by our definition.

We omit the verifications that this isomorphism between
the bifunctors $\mathbf E_{\mathbf V}(Q,M)$ and
$\Ext(0\semitimes Q,{_0M})$ of $M$ and $Q$ is natural in
both arguments, and that it preserves the abelian group
operations. \end{proof}

\section{Clone Cohomology}\label{S:clone}

The partial cochain complex which we introduced in
\S\ref{S:homologyobject}, leading to the cohomology group
we have identified as the set of extensions $\mathbf
E_{\mathbf V}(A,M)$, can be extended to a full positive
cochain complex (i.e., with objects $C^i\in\Ab[A,\mathbf
V]$ for all $i\geq 0$) and we will do so in this section.
As we do, we want to make some small changes in our
development. We want to develop the theory in such a way
that, as we did in defining $\mathbf E_{\mathbf V}(Q,M)$
in \S\ref{S:bifunctor}, the resulting cohomology objects
are bifunctors contravariant in $\Ov[A,\mathbf V]$ and
covariant in
$\Ab[A,\mathbf V]$. In addition, we want to express the
cochain complex as a result of applying a hom functor to
a chain complex we will define. Finally, we want to make
explicit the simplicial aspect of the definition. We will
call the resulting cohomology theory \emph{clone
cohomology}, because of the role played in the definition
by the clone of the variety $\mathbf V$.

\subsection{The $A$-sets $X^i_{\mathbf V}(Q)$}As a first
step,
we will define functors $X^i_{\mathbf
V}:\Ov[A,\mathbf V]\to A\hbox{-}\Set$. The $A$-sets
$X^i_{\mathbf V}(A)$, for $i=0$, $1$, and $2$, can be seen
as the special cases of $X^i_{\mathbf V}(Q)$ where $Q$ is
the $A$-set $\lsil A,1_A\rsil$.

If $Q$ is an object of $\Ov[A,\mathbf V]$, we will define
$X^0_{\mathbf V}(Q)$ to be the underlying $A$-set of
$\lsil A\semitimes Q,\pi_Q\rsil$, with the element $\pair
aq\in{_aX^0_{\mathbf V}(Q)}$ being written as $[q]_a$.

We define $X^1_{\mathbf V}(Q)$ to be the $A$-set given by
letting $_aX^1_{\mathbf V}(Q)$ be the set of triples,
written $[v;\mathbf q]_{\mathbf a}$, where
$v\in\Clo_n\mathbf V$ for some $n$, $\mathbf a\in A^n$,
and $\mathbf q\in {_{\mathbf a}Q^{\boxtimes n}}$, and such that
$v(\mathbf a)=a$.

If $i>1$, we will define $_aX^i_{\mathbf V}(Q)$ to be the
set of $(i+2)$-tuples, written
\[[v_0,\mathbf v_1,\ldots,\mathbf v_{i-1};
\mathbf q]_{\mathbf a},\]
where $v_0$ is an element of $\Clo_{n_0}\mathbf V$,
$\mathbf v_j$ for $0<j<i$ is an $n_{j-1}$-tuple of elements
of $\Clo_{n_j}(\mathbf V)$, $\mathbf a\in A^{n_{i-1}}$, and
$\mathbf q\in{_{\mathbf a}Q^{\boxtimes n_i}}$, all for some
natural numbers $n_0$, $\ldots$, $n_{i-1}$, and such that
$v_0\mathbf v_1\ldots\mathbf v_{i-1}(\mathbf a)=a$.

Note that $v_0$, $\mathbf v_1$, $\ldots$, $\mathbf v_{i-1}$
can be considered as composable arrows of the theory
constructed from $\Clo\mathbf V$.

\subsection{$\partial$}
First, we will define $\partial$
directly, making
$A$-$\Set(X^\bullet_{\mathbf V}(Q),M)$ into a chain complex.
Given an $A$-function $\factset:X^i_{\mathbf V}(Q)\to M$, we
define $\partial\factset$ by
\[(\partial\factset)[v;\mathbf q]_{\mathbf a}=v^M_{\mathbf
a}\factset[\mathbf q]_{\mathbf a}-\factset[v^Q_{\mathbf
a}(\mathbf q)]_{v(\mathbf a)}\] if $i=0$, where $[\mathbf
q]_{\mathbf a}$ stands for
$\langle [q_1]_{a_1},\ldots,[q_n]_{a_n}\rangle$, and
otherwise, by
\begin{align*}
(\partial
\factset)[v_0,\mathbf v_1,\ldots,\mathbf v_{i-1};\mathbf q\,]_{\mathbf a}
&=(v_0)^M_{\mathbf v_1\ldots\mathbf v_{i-1}(\mathbf a)}
\factset[\mathbf v_1,\ldots,\mathbf v_{i-1};
\mathbf q]_{\mathbf a}\\
&\quad-\factset[v_0\mathbf v_1,\mathbf v_2,
\ldots,\mathbf v_{i-1}; \mathbf q\,]_{\mathbf a}\\
&\quad+\factset[v_0,\mathbf v_1\mathbf v_2,
\ldots,\mathbf v_{i-1}; \mathbf q\,]_{\mathbf a}\\
&\quad\ldots\\
&\quad+(-1)^{j}\factset[v_0,\mathbf v_1,\ldots,
\mathbf v_{j-1}\mathbf v_j,\ldots,\mathbf v_{i-1};
\mathbf q\,]_{\mathbf a}\\
&\quad\ldots\\
&\quad+(-1)^i\factset[v_0,\mathbf v_1,\ldots,\mathbf v_{i-2};
(\mathbf v_{i-1})^Q_{\mathbf a}(\mathbf q\,)]_{\mathbf v_{i-1}
(\mathbf a)},
\end{align*}
where $\factset[\mathbf v_1,\ldots,\mathbf v_{i-1};\mathbf q]_{\mathbf a}$ stands for $\langle\factset[v_{11},\mathbf v_2,\ldots,\mathbf v_{i-1};\mathbf q]_{\mathbf a},\ldots,
\factset[v_{1n_0},\mathbf v_2,\ldots,\mathbf v_{i-1};\mathbf q]_{\mathbf a}\rangle$.

A straightforward computation
shows that $\partial\partial=0$, making
$A$-$\Set(X^\bullet_{\mathbf V}(Q),M)$ into a cochain
complex, which we denote by $C^\bullet_{\mathbf V}(Q,M)$..

Another way to arrive at the same cochain complex is to
apply the free functor from $A\hbox{-}\Set$ to
$\Ab[A,\mathbf V]$ to the $A$-sets $X^i_{\mathbf V}(Q)$,
yielding objects $C_i(Q,\mathbf V)\in\Ab[A,\mathbf V]$. We
then consider the simplicial complex of objects of
$\Ab[A,\mathbf V]$ given by the $C_i(Q,\mathbf V)$ and face
maps $\partial^i_j$, for $i\geq 0$ and $0\leq j\leq i$,
defined on generators by \[\partial^0_0[q]_a=0,\] by
\[\partial^1_0[v;\mathbf q\,]_{\mathbf a}=v^M_{\mathbf
a}[\mathbf q]_{\vec a}\hbox{\ \ and\ \ }
\partial^1_1[v;\mathbf q]_{\mathbf a}=[v^Q_{\mathbf
a}(\mathbf q)]_{v(\mathbf a)},\] and by
\begin{align*}
\partial^i_j[v_0,\mathbf v_1&,\ldots,
\mathbf v_{i-1};\mathbf q]_{\mathbf a}\\
&=
\begin{cases}
(v_0)^M_{\mathbf v_1\ldots\mathbf v_{i-1}(\mathbf a)}
[\mathbf v_1,\ldots,\mathbf v_{i-1};\mathbf q]_{\mathbf
a},
&\text{for $j=0$,}\\
{}[v_0,\mathbf v_1,\ldots,\mathbf v_{j-1}\mathbf v_j,
\ldots,\mathbf v_{i-1};\mathbf q]_{\mathbf a},
&\text{for $0<j<i$, and}\\
{}[v_0,\mathbf v_1,\ldots,\mathbf v_{i-2};
(\mathbf v_{i-1})^Q_{\mathbf a}(\mathbf q\,)]_{\mathbf
v_{i-1}(\mathbf a)},
&\text{for $j=i$,}
\end{cases}
\end{align*}
for $i\geq 2$. (Edge maps can also be defined, using
projection elements of $\Clo\mathbf V$, but we have no
need to do so.) We then form $\partial$ as the alternating
sum of the face maps $\partial^i_j$ over $j$, and obtain a
chain complex $C_\bullet(Q,\mathbf V)$ such that
$\Ab[A,\mathbf V](C_\bullet(Q,\mathbf V),M)=C^\bullet_{\mathbf V}(Q,M)$, the same cochain complex we described previously.

\subsection{Clone cohomology}We define the \emph{
clone cohomology objects for $Q$, with coefficients in
$M$}, to be the cohomology groups of the cochain complex
$C^\bullet_{\mathbf V}(Q,M)$, and denote
them by $H^i_{\mathbf V}(Q,M)$. Clearly, they are bifunctors
$H^i:\Ov[A,\mathbf V]^{\text{op}}\times\Ab[A,\mathbf V]\to\Ab$.

Note that the clone cohomology objects are
defined whether or not the variety $\mathbf V$ is
congruence-modular.

\subsection{Factor sets in terms of $Q$}We previously defined factor sets of extensions of $A$ by
$M$, and then defined extensions of an $A$-overalgebra
$Q$ by $M$. We have not defined factor sets in terms of
$Q$ yet, and will not do so in detail. The key fact which
allows us to relate our previous work with factor sets,
and the objects $X^\bullet_{\mathbf V}(Q)$ and
$C^\bullet_{\mathbf V}(Q,M)$, is as follows:

\begin{theorem} For $i=0$, $1$, and $2$, the three
functors to $\Ab[A,\mathbf V]$,
\[A\text{-}\Set(X^i_{\mathbf V}(A\semitimes
Q),{_{\pi_Q}{\Res M}}),\]
\[A\text{-}\Set(X^i_{\mathbf V}(Q),M),\text{ and}\]
\[C^i_{\mathbf V}(Q,M)\]
are naturally isomorphic as bifunctors in $Q$ and $M$.
\end{theorem}

\begin{proof}
The first and second functors are naturally isomorphic,
because of an adjunction between the functor of
restriction of $A$-sets along $\pi_Q$ and the functor from
$(A\semitimes Q)$-sets to $A$-sets given by sending an
$(A\semitimes Q)$-set $\lsil B,\pi\rsil$ to the $A$-set
$\lsil B,\pi_Q\circ\pi\rsil$.

The second and third functors are naturally isomorphic,
because of the adjunction between the free functor from
$A$-$\Set$ to $\Ab[A,\mathbf V]$, and the corresponding
forgetful functor.
\end{proof}

\section{Relative Clone Cohomology
$\mathbf V$}\label{S:relative}

The clone cohomology functors $H^i_{\mathbf V}(Q,M)$ are
specific to a given variety of algebras $\mathbf V$,
such that $Q$ and $M$ are totally in $\mathbf V$.
However, the construction of the chain complex used in
computing them can use a smaller variety $\mathbf V'$, if
$M$ is totally in $\mathbf V'$. That is, in that case, we
can form the free object in $\Ab[A,\mathbf V']$ on
$A$-set of generators $X^i_{\mathbf V}(Q)$, which we
denote by $C_{i,\mathbf V'}(Q,\mathbf V)$, rather than the free
object in $\Ab[A,\mathbf V]$ which we denoted by
$C_i(Q,\mathbf V)$. Because $M$ is totally in $\mathbf V'$, the
hom functor $\Ab[A,\mathbf V](-,M)$ takes these objects
to the same abelian group.  The
definition of $\partial$ also makes sense. As a result,

\begin{theorem} If $M$ is totally in $\mathbf V'$,
then the objects $C^i_{\mathbf V}(Q,M)$ and $H^i_{\mathbf V}(Q,M)$ are also totally in $\mathbf V'$. \end{theorem}

If $Q$ is also totally in $\mathbf V'$, then the definition
of $\partial$ for the chain complex $C_\bullet(Q,\mathbf V')$ makes sense, and there is an obvious homomorphism
$\pi_i$ of $C_i(Q,\mathbf V)$ onto $C_i(Q,\mathbf V')$ for each
$i$. As a result, we have
a short exact sequence of complexes in $\Ab[A,\mathbf V]$,
\[\CD K_0(Q,\mathbf V',\mathbf V) @<<< K_1
(Q,\mathbf V',\mathbf V) @<<< K_2(Q,\mathbf V',\mathbf V) @<<<
\ldots \\
@VVV @VVV @VVV \\
C_{0,\mathbf V'}(Q,\mathbf V) @<<< C_{1,\mathbf V'}(Q,\mathbf V)
@<<< C_{2,\mathbf V'}(Q,\mathbf V) @<<< \ldots \\
@VV\pi_0V @VV\pi_1V @VV\pi_2V \\
C_0(Q,\mathbf V') @<<< C_1(Q,\mathbf V') @<<< C_2(Q,\mathbf V')
@<<< \ldots \endCD
\]
where the $K_i(Q,\mathbf V',\mathbf V)$ are the kernels of
the onto homomorphisms $\pi_i$. Note that the vertical
short exact sequences split. For, each generator of
$C_i(Q,\mathbf V')$ can be mapped to a preimage in
$C_{i,\mathbf V'}(Q,\mathbf V)$, leading to a splitting of
$\pi_i$.

Applying the hom functor $\Ab[A,\mathbf V](-,M)$ and
denoting $\Ab[A,\mathbf V](K_i(Q,\mathbf V',\mathbf V),M)$
by $C^i_{\mathbf V',\mathbf V}(Q,M)$, we obtain a short exact
sequence of complexes of objects of $\Ab[A,\mathbf V']$:
\[\CD C^0_{\mathbf V',\mathbf V}(Q,M) @>>> C^1_{\mathbf
V',\mathbf V}(Q,M)
 @>>> C^2_{\mathbf V',\mathbf V}(Q,M) @>>> \ldots \\ @AAA @AAA
@AAA \\ C^0_{\mathbf V}(Q,M) @>>> C^1_{\mathbf V}(Q,M) @>>>
C^2_{\mathbf V}(Q,M) @>>> \ldots \\
@AAA @AAA @AAA \\
C^0_{\mathbf V'}(Q,M) @>>> C^1_{\mathbf V'}(Q,M) @>>>
C^2_{\mathbf V'}(Q,M) @>>> \ldots \endCD \]
where, again, the vertical sequences are exact.

We define the \emph{relative clone cohomology of $Q$, with
coefficients in $M$, with respect to the inclusion $\mathbf V'\subseteq\mathbf V$}, to be the cohomology objects of the
complex $C^\bullet_{\mathbf V',\mathbf V}(Q,M)$, and denote
these objects by $H^i_{\mathbf V',\mathbf V}(Q,M)$. It is
clear that $C^0_{\mathbf V}(Q,M)$ and $C^0_{\mathbf V'}(Q,M)$
are isomorphic; thus, $C^0_{\mathbf V',\mathbf V}(Q,M)=H^0_{\mathbf V',\mathbf V}(Q,M)=0$.  We then see that
there is a long exact sequence of cohomology groups
\[0\to H^1_{\mathbf V'}(Q,M)\to H^1_{\mathbf V}(Q,M)\to
H^1_{\mathbf V',\mathbf V}(Q,M)\to H^2_{\mathbf
V'}(Q,M)\to\ldots\] relating the three sets of cohomology
objects
$H^\bullet_{\mathbf V'}(Q,M)$, $H^\bullet_{\mathbf V}(Q,M)$,
and $H^\bullet_{\mathbf V',\mathbf V}(Q,M)$.

If we have three varieties $\mathbf V''\subseteq\mathbf V'\subseteq\mathbf V$ such that $Q$ and $M$ are totally in
$\mathbf V''$, then similar methods, and standard methods
from homological algebra, yield a long exact sequence
\[0\to H^1_{\mathbf V'',\mathbf V'}(Q,M)\to H^1_{\mathbf
V'',\mathbf V}(Q,M)\to H^1_{\mathbf V',\mathbf V}(Q,M)\to
H^2_{\mathbf V'',\mathbf V'}(Q,M)\to\ldots\] relating the
relative cohomology objects.

 \section*{Discussion}
The study of abelian extensions, and
the recognition that they form a cohomology group, has a
considerable history, as does the investigation of
cohomology theories in general. We have focussed our
attention on the abelian group of extensions, and have
defined a cohomology theory that has this algebra as
cohomology group in dimension one.  Many questions remain
to be answered about this new cohomology theory, and about
its relationship to previously-studied theories of
cohomology of algebras. We will raise some of those
questions in this section.

Categorical algebraists have invented a cohomology theory
called comonadic cohomology. The theory can be applied to
any comonad, but usually, the comonad being used when
this theory is mentioned is a comonad derived from
the free and forgetful functors for the variety in
question.  See
\cite{1}, \cite{3}, and \cite{2}. As in clone
cohomology, the theory gives the cohomology group of an
$A$-overalgebra (actually, an algebra over $A$, as the
theory is usually developed) totally in a variety $\mathbf V$
to which $A$ belongs, with coefficients in an $A$-module
$M$ totally in $\mathbf V$ (or, as usually expressed, in a
Beck module over $A$). As in clone cohomology, the
resolution that gives rise to the cohomology objects comes
from taking alternating sums of face maps of a simplicial
complex.

We are led to ask, what is the relationship of the clone
cohomology objects to the comonadic
cohomology objects? 
 In \cite{3}, the comonadic
cohomology groups in dimension 0 and 1 were studied and
interpreted. This was done in the generality of $\mathbf V$
an arbitrary (i.e., not necessarily congruence-modular)
variety of algebras of some type. The interpretation of
the group in dimension 1 is different from our
interpretation of the clone cohomology group in dimension
1, but, we have proved directly that, for $\mathbf V$ a
congruence-modular variety of algebras, the groups are
naturally isomorphic. We have not discussed dimension 0
in this paper, but the groups in dimension 0 are also
isomorphic. (We have not included the proofs of these
results in this paper.)

This raises the question, are the groups isomorphic in
all dimensions, when $\mathbf V$ is congruence-modular?
More generally, are they isomorphic when $\mathbf V$ is not
congruence-modular? It is not hard to show that they are
isomorphic in dimension 0, but the answer is not known
for higher dimensions, or even for dimensions higher than
1 in the congruence-modular case.

For clone cohomology, and for comonadic cohomology
for that matter, there is the question, when are the
resolutions used in the derivation of the cohomology
groups exact? Note that in both cases, the resolutions are
free. Thus, when the resolutions are exact, the cohomology
functors will be derived functors in the standard sense.
Otherwise, for comonadic cohomology, there is at least a
uniqueness theorem (\cite{2}) which characterizes
the  cohomology functors. For clone cohomology, there is
not yet such a theorem.

For some well-known cases, such as the variety of
groups, the comonadic and clone cohomology groups
coincide with the standard cohomology groups which can be
defined as derived functors in the usual sense. Thus, we
can say that we have characterized these groups using two
different universal properties. What is the significance of
the fact that the two different universal properties arrive
at the same answer?

Finally, we should mention that there is a formal, at
least, resemblence between the derivation of the clone
cohomology objects, and the construction of the
cohomology groups of a category. What is the significance
of this
resemblence?

To our knowledge, none of the questions have yet been
answered.

\subsection{Acknowledgement} It is a pleasure to
acknowledge helpful discussions and correspondence
about this topic with Saunders Mac Lane.

\end{document}